\documentstyle[12pt,amssymb]{amsart}

\numberwithin{equation}{section}

\newtheorem{lem}{Lemma}[section]
\newtheorem{prp}[lem]{Proposition}
\newtheorem{thm}[lem]{Theorem}

\begin{document}

\title[hypergeometric polynomials in several variables]{Properties of
some families of hypergeometric orthogonal polynomials in several variables}
\author[J. F. van Diejen]{J. F. van Diejen}
\address{Department of Mathematical Sciences, The University of Tokyo,
3-8-1 Komaba, Meguro-ku, Tokyo~153, Japan}


\thanks{Work supported by the Japan Society for the Promotion of Science
(JSPS) and by a Monbusho Grant-in-Aid.}

\date{(March 1996)}
\maketitle

\begin{abstract}
Limiting cases are studied of the Koornwinder-Macdonald multivariable
generalization of the Askey-Wilson polynomials. We recover
recently and not so recently introduced families of hypergeometric orthogonal
polynomials in several variables consisting of multivariable
Wilson, continuous Hahn and Jacobi type polynomials, respectively.
For each class of polynomials we provide systems of difference
(or differential) equations, recurrence relations, and expressions for
the norms of the polynomials in terms of the norm of the constant polynomial.
\end{abstract}

%\newpage
\section{Introduction}\label{sec1}
It is to date over a decade ago that Askey and Wilson released
their famous memoir \cite{ask-wil:some}, in which they introduced a
four-parameter family of basic hypergeometric polynomials nowadays
commonly referred to as the Askey-Wilson polynomials
\cite{gas-ram:basic}. These polynomials, which are defined explicitly
in terms of a terminating ${}_4\phi_3$ series,
have been shown to exhibit a number of interesting properties.
Among other things, it was demonstrated that they
satisfy a second order difference equation, a three-term recurrence
relation, and that---in a suitable parameter regime---they constitute
an orthogonal system with respect to an explicitly given positive weight
function with support on a finite interval
(or on the unit circle, depending on how
the coordinates are chosen).

Many (basic) hypergeometric orthogonal polynomials studied in the
literature arise as special (limiting) cases of the Askey-Wilson
polynomials and have been collected in
the so-called ($q$-)Askey scheme \cite{ask-wil:some,koe-swa:askey-scheme}.
For instance, if the step size parameter of the difference equation is
scaled to zero, then the Askey-Wilson polynomials go over in Jacobi
polynomials: well-known classical hypergeometric orthogonal
polynomials satisfying a second order differential equation
instead of a difference equation. One may
also consider the transition from orthogonal polynomials on a
finite interval to orthogonal polynomials on a (semi-)infinite
interval. This way
one arrives at Wilson polynomials (semi-infinite interval)
and at continuous Hahn
polynomials (infinite interval).

The purpose of the present paper is to generalize this state of
affairs from one to several variables. Starting point is a recently
introduced multivariable generalization of the Askey-Wilson
polynomials, found for special parameters by Macdonald
\cite{mac:orthogonal} and in full generality (involving five parameters)
by Koornwinder \cite{koo:askey}.
By means of limiting transitions similar to those
in the one-variable case, we arrive at multivariable Jacobi polynomials
\cite{vre:formulas,deb:systeme}
(see also \cite{bee-opd:certain} and reference therein) and at
multivariable Wilson and continuous Hahn polynomials
\cite{die:multivariable}.

The ($q$-)Askey scheme involves many more limits and special cases of the
Askey-Wilson polynomials than those described above. For instance, one
also considers transitions from certain polynomials in the scheme to
similar polynomials with less parameters and transitions from
polynomials with a continuous orthogonality measure to polynomials
with a discrete orthogonality measure. Such transitions (or rather
their multivariable analogues) will not be considered here. We refer
instead to \cite[Sec. 5.2]{die:commuting} for the transition from
multivariable Askey-Wilson polynomials to Macdonald's $q$-Jack
polynomials (i.e., multivariable $q$-ultraspherical
polynomials) \cite{mac:symmetric} (as an example of a limit leading to
similar polynomials but with less parameters), and to \cite{sto-koo:limit}
for the transition from multivariable Askey-Wilson polynomials to
multivariable big and little $q$-Jacobi polynomials
\cite{sto:multivariable} (as an example of a limit leading to
multivariable polynomials with a discrete orthogonality measure).

Whenever one is dealing with orthogonal polynomials an important
question arises as to the explicit computation of the normalization constants
converting the polynomials into an orthonormal system. For Jacobi
polynomials calculating the orthonormalization constants boils down to
the evaluation of (standard) beta integrals, whereas Askey-Wilson
polynomials give rise to $q$-beta integrals. In the case of several
variables one has to deal with Selberg type integrals (Jacobi case)
and $q$-Selberg type integrals (Askey-Wilson case), respectively. For
these multiple integrals explicit evaluations have been conjectured by
Macdonald that were
recently checked using techniques involving so-called shift operators
\cite{opd:some,hec-sch:harmonic,che:double,mac:affine}.
(Roughly speaking these shift operators allow
one to relate the values of the ($q$-)Selberg integral for different
values of the parameters separated by unit shifts;
the integral can then be solved, first for nonnegative integer-valued
parameters by shifting the parameters
to zero in which case the integrand becomes trivial,
and then for arbitrary nonnegative parameters using an analyticity argument
(viz. Carlson's theorem).)

Very recently, the author observed that Koornwinder's second order
difference equation for the multivariable Askey-Wilson polynomials may
be extended to a system of $n$ ($=$ number of variables) independent
difference equations \cite{die:commuting} and that the polynomials
also satisfy a system of $n$ independent recurrence relations
\cite{die:self-dual}. (To date a complete proof for these recurrence
relations is only available for a self-dual four-parameter subfamily of the
five-parameter multivariable Askey-Wilson polynomials.)
It turns out that the recurrence relations, combined with the known
evaluation for the norm of the unit polynomial (i.e., the constant term
integral) \cite{gus:generalization,kad:proof}, may also be used to verify
Macdonald's formulas for the orthonormalization constants of the
multivariable Askey-Wilson polynomials \cite{die:self-dual}.
Below, we will use these results to arrive
at systems of difference (or differential) equations, recurrence
relations and expressions for the orthonormalization constants, for
all three limiting cases
of the multivariable Askey-Wilson polynomials
considered in this paper (Wilson, continuous Hahn and Jacobi type).

We would like to emphasize that much of the presented material
admits a physical interpretation in terms
of Calogero-Sutherland type exactly solvable
quantum $n$-particle models related to classical root systems
\cite{ols-per:quantum} or their Ruijsenaars type difference versions
\cite{rui:complete,rui:finite,die:difference}.
The point is that the second order
differential equation for the multivariable Jacobi polynomials
may be seen as the eigenvalue equation for
a trigonometric quantum Calogero-Sutherland system
related to the root system $BC_n$ \cite{ols-per:quantum}.
{}From this viewpoint the second order
difference equation for the multivariable Askey-Wilson polynomials
corresponds to the eigenvalue equation for a Ruijsenaars type
difference version of the $BC_n$-type quantum
Calogero-Sutherland system \cite{die:diagonalization}.
The transitions to the multivariable continuous Hahn
and Wilson polynomials amount to rational limits leading to
(the eigenfunctions of) similar difference
versions of the $A_{n-1}$-type rational Calogero model with harmonic term
(continuous Hahn case) and its $B(C)_n$-type counterpart (Wilson case)
\cite{die:multivariable}.
For further details regarding these connections with the Calogero-Sutherland
and Ruijsenaars type quantum integrable
$n$-particle systems the reader is referred to
\cite{die:difference,die:multivariable,die:diagonalization}.

The material is organized as follows.
First we define our families of multivariable (basic) hypergeometric
polynomials in Section~\ref{sec2} and recall their second
order difference equation (Askey-Wilson, Wilson, continuous Hahn type)
or second order differential equation (Jacobi type)
in Section~\ref{sec3}. Next, in Section~\ref{sec4}, limit transitions
from the Askey-Wilson type family to the Wilson, continuous Hahn
and Jacobi type families are discussed.
We study the behavior of our recently introduced systems of
difference equations and recurrence
relations for the multivariable Askey-Wilson type polynomials
with respect to
these limits in Sections~\ref{sec5} and \ref{sec6}, respectively.
The recurrence relations for the Wilson, continuous Hahn
and Jacobi type polynomials thus obtained in Section~\ref{sec6} are then
employed in Section~\ref{sec7}
to derive explicit expressions for the (squared) norms of the
corresponding polynomials
in terms of the (squared) norm of the unit polynomial.

\vfill

\section{Multivariable (basic) hypergeometric polynomials}\label{sec2}
In this section multivariable versions of some orthogonal
families of (basic) hypergeometric polynomials are characterized.
The general idea of the construction (which is standard, see
e.g. \cite{vre:formulas,mac:orthogonal,koo:askey,sto-koo:limit})
is to start with an algebra
of (symmetric) polynomials ${\cal H}$ spanned by a basis of
(symmetric) monomials
$\{ m_\lambda \}_{\lambda \in \Lambda}$, with the set $\Lambda$
labeling the basis elements being partially ordered in such a way that
for all $\lambda\in\Lambda$ the subspaces
${\cal H}_\lambda \equiv \text{Span}\{ m_\mu \}_{\mu\in\Lambda ,\mu
\leq\lambda}$ are finite-dimensional.
It is furthermore assumed that the space ${\cal H}$ is endowed with an $L^2$
inner product $\langle \cdot ,\cdot \rangle_{\Delta}$ characterized by
a certain weight function $\Delta$.
To such a configuration we associate a
basis $\{p_\lambda\}_{\lambda\in\Lambda}$ of ${\cal H}$ consisting of
the polynomials $p_\lambda$, $\lambda\in\Lambda$, determined
(uniquely) by the two conditions
\begin{itemize}
\item[{\em i.}] $\displaystyle p_\lambda  = m_\lambda  +
\sum_{\mu \in \Lambda, \mu < \lambda }\;
 c_{\lambda ,\mu }\: m_{\mu}$,
\ \ \ \ \ $\displaystyle c_{\lambda ,\mu}\in {\Bbb C}$;
\item[{\em ii.}]
$\displaystyle  \langle p_\lambda  ,  m_{\mu}\rangle_\Delta  =0$
\ \ if \ \  $\displaystyle \mu < \lambda$.
\end{itemize}
In other words, the polynomial $p_\lambda$ consists of the monomial
$m_\lambda$ minus its orthogonal
projection with respect to the
inner product $\langle \cdot ,\cdot \rangle_{\Delta}$
onto the finite-dimensional subspace
$\text{Span}\{ m_\mu \}_{\mu\in\Lambda ,\mu <\lambda}$.
By varying the concrete choices for the space ${\cal H}$, the basis
$\{ m_\lambda \}_{\lambda \in \Lambda}$ and the inner product
$\langle \cdot ,\cdot \rangle_{\Delta}$, we recover certain
(previously introduced) multivariable
generalizations of the Askey-Wilson, Wilson, continuous Hahn and
Jacobi polynomials, respectively.
Below we will specify the relevant
data determining these families.
The fact that in the case of one variable the corresponding polynomials
$p_\lambda$ indeed reduce to the well-known one-variable
polynomials studied extensively in the literature is immediate from the weight
function. The normalization for the polynomials is determined by the
fact that (by definition)
$p_\lambda$ is monic in the sense that the coefficient of
the leading monomial $m_\lambda$ in $p_\lambda$ is equal to one.

It turns out that in all of our cases the basis $\{ m_\lambda\}_{\lambda \in
\Lambda}$ can be conveniently expressed in terms of the monomial symmetric
functions
\begin{equation}\label{monsym}
m_{sym,\lambda}(z_1,\ldots ,z_n)=
\sum_{\mu \in S_n (\lambda )} z_1^{\mu_1}\cdots z_n^{\mu_n},\;\;\;\;\;
\lambda \in \Lambda ,
\end{equation}
where
\begin{equation}\label{cone}
\Lambda  = \{ \lambda \in {\Bbb Z}^n\; |\;
\lambda_1\geq\lambda_2\geq \cdots \geq \lambda_n\geq 0 \} .
\end{equation}
In \eqref{monsym} the summation is meant
over the orbit of $\lambda$ under the action of
the permutation group $S_n$ (acting on the vector components
$\lambda_1,\ldots ,\lambda_n$).
As partial order of the integral cone $\Lambda$ \eqref{cone} we will always
take the dominance order defined by
\begin{equation}\label{po}
\mu \leq \lambda \;\;\;\; \text{iff}\;\;\;\;
\sum_{j=1}^m \mu_j \leq \sum_{j=1}^m \lambda_j \;\; \text{for}\;\;
m=1,\ldots ,n
\end{equation}
(and $\mu <\lambda$ iff $\mu \leq \lambda$ and $\mu \neq \lambda$).

{\em Note:} In order to avoid confusion between the various
families we will often
equip the polynomials and other objects
of interest with the superscripts `AW', `W', `cH' or `J' to indicate
Askey-Wilson, Wilson, continuous Hahn or Jacobi type polynomials,
respectively. Sometimes, however, these superscripts will be suppressed when
discussing more general properties of the polynomials
that hold simultaneously for all families.

\subsection{Askey-Wilson type}\label{AWdef}
To arrive at multivariable Askey-Wilson type polynomials one considers a
space ${\cal H}^{AW}$ consisting of even and permutation invariant
trigonometric polynomials. Specifically, the space ${\cal H}^{AW}$
is spanned by the monomials
\begin{equation}\label{AWmonom}
m_\lambda^{AW} (x) = m_{sym,\lambda}
(e^{i\alpha x_1}+e^{-i\alpha x_1},\cdots ,
e^{i\alpha x_n}+e^{-i\alpha x_n}),\;\;\;\;\; \lambda\in\Lambda
\end{equation}
(with $\Lambda$ given by \eqref{cone}).
The relevant inner product on ${\cal H}^{AW}$ is determined by
\begin{equation}\label{ip1}
\langle m^{AW}_\lambda ,m^{AW}_{\mu} \rangle_{\Delta^{AW}} =
\int_{-\pi/\alpha}^{\pi/\alpha}\!\!\!\!\!\!\cdots
               \int_{-\pi/\alpha}^{\pi/\alpha}
m^{AW}_\lambda (x)\, \overline{m^{AW}_{\mu}(x)}\,  \Delta^{AW} (x)\,
dx_1\cdots dx_n ,
\end{equation}
with the weight function reading
\begin{eqnarray}\label{weightAW}
\Delta^{AW} (x) &=&
\prod_{\stackrel{1\leq j< k \leq n}
                {\varepsilon_1 ,\varepsilon_2 =\pm 1}}
\frac{(e^{i\alpha (\varepsilon_1 x_j
           +\varepsilon_2 x_k)}; q)_{\infty}}
     {(t\, e^{i\alpha (\varepsilon_1 x_j
           +\varepsilon_2 x_{k})}; q)_{\infty}} \\
& &  \times \prod_{\stackrel{1\leq j\leq n}{\varepsilon =\pm 1}}
\frac{( e^{2i\alpha \varepsilon x_j}; q)_\infty    }
     {(t_0\, e^{i\alpha \varepsilon x_j},\,
       t_1\, e^{i\alpha \varepsilon x_j},\,
       t_2\, e^{i\alpha \varepsilon x_j},\,
       t_3\, e^{i\alpha \varepsilon x_j};\, q)_\infty    } .  \nonumber
\end{eqnarray}
Here $(a;q)_\infty \equiv \prod_{m=0}^\infty (1-aq^m)$,
$(a_1,\ldots ,a_r;q)_\infty\equiv (a_1;q)_\infty\cdots (a_r;q)_\infty$
and the parameters are assumed to satisfy the constraints
\begin{equation}
\alpha > 0, \;\;\;\; 0< q<1,\;\;\;\; -1\leq t\leq 1, \;\; \;\;
|t_r|\leq 1\;\; (r=0,1,2 ,3),
\end{equation}
with possible
non-real parameters $t_r$ occurring in complex conjugate pairs
and pairwise products of the $t_r$ being $\neq 1$.
For the weight function in \eqref{weightAW}
the polynomials $p_\lambda$ determined by the
Conditions {\em i.} and {\em ii.} (above) were introduced by
Macdonald \cite{mac:orthogonal} (for special parameters) and Koornwinder
\cite{koo:askey} (for general parameters).
In the special case of one variable ($n=1$)
these polynomials reduce to monic Askey-Wilson polynomials
\cite{ask-wil:some,koe-swa:askey-scheme}
\begin{equation}\label{bhrepAW}
p^{AW}_l(x) =
\frac{(t_0 t_1, t_0 t_2, t_0 t_3; q)_l}{t_0^l (t_0 t_1 t_2 t_3
q^{l-1};q)_l }\;
{}_4\phi_3
\left(
\begin{matrix}
q^{-l},\; t_0 t_1 t_2 t_3 q^{l-1},\;
t_0 e^{i\alpha x},\; t_0 e^{-i\alpha x} \\ [0.5ex]
t_0 t_1 ,\; t_0 t_2 ,\; t_0 t_3
\end{matrix} \; ; q,q \right) .
\end{equation}

\subsection{Wilson type}\label{Wdef}
In the Wilson case the appropriate space ${\cal H}^W$ consists of
even and permutation invariant polynomials and is
spanned by the monomials
\begin{equation}\label{Wmonom}
m^W_\lambda (x) = m_{sym,\lambda}
(x_1^2,\cdots ,x_n^2),\;\;\;\;\; \lambda\in\Lambda .
\end{equation}
The inner product on ${\cal H}^W$ is now determined by
\begin{equation}\label{ip2}
\langle m^W_\lambda ,m^W_{\mu} \rangle_{\Delta^W} =
\int_{-\infty}^{\infty}\!\!\!\!\!\!\cdots
               \int_{-\infty}^{\infty}
m^W_\lambda (x)\, \overline{m^W_{\mu}(x)}\,  \Delta^W (x)\,
dx_1\cdots dx_n ,
\end{equation}
with the weight function taken to be of the form
\begin{eqnarray}\label{weightW}
\Delta^W (x) &=&
\prod_{\stackrel{1\leq j< k \leq n}
                {\varepsilon_1 ,\varepsilon_2 =\pm 1}}
\frac{\Gamma (\nu + i(\varepsilon_1 x_j
           +\varepsilon_2 x_{k}) )}
     {\Gamma (i(\varepsilon_1 x_j
           +\varepsilon_2 x_{k}) )} \\
& &  \times \prod_{\stackrel{1\leq j\leq n}{\varepsilon =\pm 1}}
\frac{\Gamma (\nu_0+i\varepsilon x_j)
\Gamma (\nu_1+i\varepsilon x_j)
\Gamma (\nu_2+i\varepsilon x_j)
\Gamma (\nu_3+i\varepsilon x_j)}{\Gamma (2i\varepsilon x_j )} .  \nonumber
\end{eqnarray}
Here $\Gamma (\cdot)$ denotes the gamma function and the parameters
are such that
\begin{equation}\label{Wcond}
\nu \geq 0, \;\;\;\;\;\; \text{Re}(\nu_r) > 0  \;\; (r=0,1,2 ,3),
\end{equation}
with possible non-real parameters $\nu_r$ occurring in complex conjugate pairs.
For this weight function (and real parameters)
the polynomials $p_\lambda$
were introduced in \cite{die:multivariable}.
In the case of one variable
they reduce to monic Wilson polynomials
\cite{ask-wil:some,koe-swa:askey-scheme}
\begin{eqnarray}\label{hrepW}
p^W_l(x) &=&
\frac{(\nu_0+\nu_1, \nu_0+\nu_2, \nu_0+\nu_3)_l}
     {(-1)^l (\nu_0+\nu_1+\nu_2+\nu_3+l-1)_l } \times\\
& & {}_4F_3
\left(
\begin{matrix}
-l,\; \nu_0+\nu_1+\nu_2+\nu_3+l-1,\;
\nu_0+ix,\; \nu_0-ix \\ [0.5ex]
\nu_0+\nu_1 ,\; \nu_0+\nu_2 ,\; \nu_0+\nu_3
\end{matrix} \; ; 1\right) . \nonumber
\end{eqnarray}

\subsection{Continuous Hahn type}\label{cHdef}
The space ${\cal H}^{cH}$ is very similar to that of the
Wilson case but instead of only the even sector it now consists of {\em all}
permutation invariant polynomials. The monomial basis for the
space ${\cal H}^{cH}$ then becomes
\begin{equation}\label{cHmonom}
m^{cH}_\lambda (x) = m_{sym,\lambda}
(x_1,\cdots ,x_n),\;\;\;\;\; \lambda\in\Lambda .
\end{equation}
The inner product is of the same form as for the Wilson case
\begin{equation}\label{ip3}
\langle m^{cH}_\lambda ,m^{cH}_{\mu} \rangle_{\Delta^{cH}} =
\int_{-\infty}^{\infty}\!\!\!\!\!\!\cdots
               \int_{-\infty}^{\infty}
m^{cH}_\lambda (x)\, \overline{m^{cH}_{\mu}(x)}\,  \Delta^{cH} (x)\,
dx_1\cdots dx_n
\end{equation}
but now with a weight function given by
\begin{eqnarray}\label{weightcH}
\Delta^{cH} (x) &=&
\prod_{1\leq j< k \leq n}
\left( \frac{\Gamma (\nu + i(x_j -x_{k}) )}
     {\Gamma (i(x_j -x_{k}) )}
\frac{\Gamma (\nu + i(x_k -x_{j}) )}
     {\Gamma (i(x_k -x_{j}) )} \right) \\
& &\times \prod_{1\leq j\leq n}
\Gamma (\nu^+_0+i x_j) \Gamma (\nu^+_1+i x_j)
\Gamma (\nu^-_0-i x_j) \Gamma (\nu^-_1-i x_j) \nonumber ,
\end{eqnarray}
where
\begin{equation}\label{cHcond}
\nu \geq 0, \;\;\;\;\;\;
\text{Re}(\nu^\pm_0), \text{Re}(\nu^\pm_1) > 0,\;\;\;\;
\nu_0^-=\overline{\nu_0^+},\;\;
\nu_1^-=\overline{\nu_1^+}.
\end{equation}
Just as in the case of Wilson type polynomials,
the polynomials $p_\lambda^{cH}$ corresponding to the weight function
\eqref{weightcH}
were introduced in \cite{die:multivariable}. For $n=1$ they reduce
to monic continuous Hahn polynomials
\cite{ask-wil:some,koe-swa:askey-scheme}
\begin{eqnarray}\label{hrepcH}
p^{cH}_l(x) &=&
 \frac{i^l\; (\nu^+_0+\nu^-_0, \nu^+_0+\nu^-_1 )_l}
     {(\nu^+_0+\nu^-_0+\nu^+_1+\nu^-_1+l-1)_l } \times \\
& &\;\;\;\;\;\;\;\;\;\; {}_3F_2
\left(
\begin{matrix}
-l,\; \nu^+_0+\nu^-_0+\nu^+_1+\nu^-_1+l-1,\;
\nu^+_0+ix \\ [0.5ex]
\nu^+_0+\nu^-_0 ,\; \nu^+_0+\nu^-_1
\end{matrix} \; ; 1\right) . \nonumber
\end{eqnarray}

\subsection{Jacobi type}\label{Jdef}
The space ${\cal H}^J$ and the basis
$\{ m^J_\lambda \}_{\lambda \in \Lambda}$ are the same as for the
Askey-Wilson type. The inner product is also of the form given there
(cf. \eqref{ip1}) but the weight function gets replaced by
\begin{eqnarray}\label{weightJ}
\Delta^J (x) &=&
\prod_{1\leq j< k \leq n}
\left|\, \sin \frac{\alpha}{2} (x_j+x_{k})
  \sin \frac{\alpha}{2}(x_j-x_{k}) \,\right|^{2\nu} \\
& &\times \prod_{1\leq j\leq n}
\left|\, \sin (\frac{\alpha}{2}x_j)\, \right|^{2\nu_0}
\left| \, \cos (\frac{\alpha}{2}x_j) \, \right|^{2\nu_1}  , \nonumber
\end{eqnarray}
with
\begin{equation}
\alpha > 0, \;\;\;\;\;\;  \nu \geq 0,\;\;\;\; \nu_0,\nu_1 > -1/2.
\end{equation}
In this case the corresponding polynomials $p_\lambda$ were first introduced
by Vretare \cite{vre:formulas} (see also
\cite{deb:systeme,bee-opd:certain}).
For $n=1$ they reduce to monic Jacobi
polynomials \cite{abr-ste:handbook,koe-swa:askey-scheme}
\begin{equation}\label{hrepJ}
p^J_l(x) =
 \frac{2^{2l}\; (\nu_0+1/2)_l}
     {(\nu_0+\nu_1+l)_l } \;\:
{}_2F_1
\left(
\begin{matrix}
-l,\; \nu_0+\nu_1+l \\ [0.5ex]
\nu_0+1/2 \end{matrix} \; ; \; \sin^2 \left(\frac{\alpha x}{2}\right) \right) .
\end{equation}

{\em Remark:}
In the explicit formulas for the polynomials when $n=1$ we have used their
standard representations in terms of terminating (basic)
hypergeometric series \cite{ask-wil:some,koe-swa:askey-scheme}
\begin{eqnarray}
{}_rF_s\left( \begin{array}{c}
             a_1,\ldots ,a_r \\ b_1,\ldots ,b_s
            \end{array} ; z \right) &=&
\sum_{k=0}^\infty
\frac{(a_1,\ldots ,a_r)_k}{(b_1,\ldots ,b_s)_k} \frac{z^k}{k!} ,\\
{}_{s+1}\phi_s\left( \begin{array}{c}
             a_1,\ldots ,a_{s+1} \\ b_1,\ldots ,b_s
            \end{array} ; q,z \right) &=&
\sum_{k=0}^\infty
\frac{(a_1,\ldots ,a_{s+1};q)_k}{(b_1,\ldots ,b_s;q)_k}
\frac{z^k}{(q;q)_k},
\end{eqnarray}
where we have used
Pochhammer symbols and $q$-shifted factorials defined by
\begin{equation*}
(a_1,\ldots ,a_r)_k= (a_1)_k\cdots (a_r)_k,\;\;\;
(a_1,\ldots ,a_r;q)_k= (a_1;q)_k\cdots (a_r;q)_k
\end{equation*}
with $(a)_0=(a;q)_0=1$ and
\begin{equation*}(a)_k=a(a+1)\cdots (a+k-1),\;\;
(a;q)_k=(1-a)(1-aq)\cdots (1-aq^{k-1})
\end{equation*}
for $k=1,2,3,\dots$.

\section{Second order difference or differential
equations}\label{sec3}
As it turns out, all families of polynomials
$\{ p_\lambda \}_{\lambda\in\Lambda}$ introduced in the previous
section satisfy an eigenvalue equation of the form
\begin{equation}\label{eveq}
D\, p_\lambda = E_\lambda \, p_\lambda, \;\;\;\;\;\; \lambda\in \Lambda,
\end{equation}
where $D:{\cal H}\rightarrow {\cal H}$ denotes a certain second order
difference operator (Askey-Wilson, Wilson and continuous Hahn case) or
a second order differential operator (Jacobi case).
Below we will list for each family
the relevant operator $D$ together with its eigenvalues
$E_\lambda$, $\lambda\in\Lambda$.
In each case the proof that the polynomials $p_\lambda$ indeed satisfy
the corresponding eigenvalue equations boils down to demonstrating
that the operator $D: {\cal H}\rightarrow {\cal H}$ maps the
finite-dimensional subspaces
${\cal H}_{\lambda} = \text{Span}\{ m_\mu \}_{\mu\in\Lambda, \mu\leq
\lambda}$ into themselves (triangularity) and that it is
symmetric with respect to the inner product
$\langle \cdot ,\cdot\rangle_{\Delta}$. In other words, one has to
show that

{\em Triangularity}
\begin{equation}\label{tria}
D\, m_\lambda = \sum_{\mu\in\Lambda ,\mu \leq\lambda}
[D]_{\lambda , \mu}\, m_\mu , \;\;\;\;\;\; \text {with}\;\;\;
[D]_{\lambda , \mu}\in {\Bbb C}
\end{equation}
and that

{\em Symmetry}
\begin{equation}\label{sym}
\langle D m_\lambda , m_\mu \rangle_\Delta =
\langle m_\lambda , D m_\mu \rangle_\Delta .
\end{equation}
It is immediate from these two properties and the definition of
the polynomial $p_\lambda$ that $D p_\lambda$ lies in ${\cal
H}_\lambda$ and is orthogonal with respect to
$\langle \cdot ,\cdot\rangle_{\Delta}$ to all monomials $m_\mu$,
$\mu\in\Lambda$ with $\mu < \lambda$.
But then comparison
with the definition of $p_\lambda$ shows that
$D p_\lambda$ must be proportional to $p_\lambda$, i.e.,
$p_\lambda$ is an eigenfunction of $D$. The corresponding eigenvalue
$E_\lambda$ is obtained via an explicit computation of the diagonal matrix
element $[D]_{\lambda ,\lambda}$ in Expansion~\eqref{tria}.

For the Jacobi case a proof of the second order differential equation
along the above lines was given by Vretare \cite{vre:formulas}.
In the Askey-Wilson case the proof was given by Macdonald
\cite{mac:orthogonal} and (in general) Koornwinder \cite{koo:askey}.
The proof for the Wilson and continuous Hahn case is very similar to
that of the Askey-Wilson case and has
been outlined in \cite{die:multivariable}.

\subsection{Askey-Wilson type}\label{sec3AW}
The second order ($q$-)difference operator
diagonalized by the polynomials $p^{AW}_\lambda$, $\lambda\in\Lambda$, is
given by
\begin{equation}\label{DAW}
D^{AW} = \sum_{1\leq j\leq n} \left(
V^{AW}_j(x) (T_{j,q}-1) \; +\; V^{AW}_{-j}(x) (T_{j,q}^{-1}-1) \right)
\end{equation}
with
\begin{eqnarray}
V^{AW}_{\pm j}(x) &=&
\frac{ \prod_{0\leq r\leq 3} (1-t_r e^{\pm i\alpha x_j})}
     { (1-e^{\pm 2i\alpha x_j})\, (1-q\, e^{\pm 2i\alpha x_j})} \\
& & \times \prod_{1\leq k \leq n,\, k\neq j}
\left( \frac{1-t e^{i\alpha (\pm x_j +x_k)}}{1-e^{i\alpha (\pm
x_j+x_k)}}\right)
\left( \frac{1-t e^{i\alpha (\pm x_j -x_k)}}{1-e^{i\alpha (\pm
x_j-x_k)}} \right) .
\nonumber
\end{eqnarray}
Here the operators $T_{j,q}$ act on trigonometric polynomials
by means of a $q$-shift of the $j$th variable
\begin{equation}\label{Tjq}
(T_{j,q}f) (e^{i\alpha x_1},\ldots ,e^{i\alpha x_n}) =
f(e^{i\alpha x_1},\ldots ,e^{i\alpha x_{j-1}},qe^{i\alpha x_j},
e^{i\alpha x_{j+1}}, \ldots ,e^{i\alpha x_n}).
\end{equation}
The eigenvalue of $D^{AW}$ on $p_\lambda^{AW}$ has the value
\begin{equation}\label{EAW}
E^{AW}_\lambda = \sum_{1\leq j\leq n} \left(
t_0t_1t_2t_3 q^{-1} t^{2n-j-1} (q^{\lambda_j}-1) +
t^{j-1}(q^{-\lambda_j}-1) \right) .
\end{equation}

\begin{prp}[\cite{koo:askey}]\label{diffAWprp}
The multivariable Askey-Wilson polynomials
$p_\lambda^{AW}$, $\lambda\in\Lambda$ \eqref{cone},
satisfy the second order difference equation
\begin{equation}\label{eveqAW}
D^{AW} p_\lambda^{AW} \; =\; E_{\lambda}^{AW}\: p_\lambda^{AW}.
\end{equation}
\end{prp}
For $n=1$, Equation \eqref{eveqAW}
reduces to the second order
difference equation for the one-variable Askey-Wilson
polynomials \cite{ask-wil:some,koe-swa:askey-scheme}
\begin{eqnarray}
\frac{ \prod_{0\leq r\leq 3} (1-t_r e^{i\alpha x})}
     { (1-e^{2i\alpha x})\, (1-qe^{2i\alpha x})}
\left( p^{AW}_l(qx)-p^{AW}_l(x) \right) &+& \\
\frac{ \prod_{0\leq r\leq 3} (1-t_r e^{-i\alpha x})}
     { (1-e^{-2i\alpha x})\, (1-qe^{-2i\alpha x})}
\left( p^{AW}_l(q^{-1} x)-p^{AW}_l(x) \right) &+& \nonumber \\
=
\left( t_0t_1t_2t_3 q^{-1} (q^{l}-1) + (q^{-l}-1)\right) p_l^{AW}(x).
& & \nonumber
\end{eqnarray}

\subsection{Wilson type}
In the case of Wilson type polynomials the difference operator
takes the form
\begin{equation}\label{DW}
D^W = \sum_{1\leq j\leq n} \left(
V^W_j(x) (T_{j}-1) \; +\; V^W_{-j}(x) (T_{j}^{-1}-1) \right)
\end{equation}
where
\begin{eqnarray}
V^W_{\pm j}(x) &=&
\frac{ \prod_{0\leq r\leq 3} (i\nu_r\pm x_j)}
     { (\pm 2ix_j )\, (\pm 2ix_j -1)} \\
& & \times \prod_{1\leq k \leq n,\, k\neq j}
\left( \frac{i\nu \pm x_j +x_k}{\pm x_j+x_k}\right)
\left( \frac{i\nu \pm x_j -x_k}{\pm x_j-x_k}\right)
\nonumber
\end{eqnarray}
and the action of $T_j$ is given by a unit shift of the $j$th variable
along the imaginary axis
\begin{equation}\label{Tj}
(T_{j}f) (x_1,\ldots ,x_n) =
f(x_1,\ldots ,x_{j-1},x_j+i,
x_{j+1}, \ldots ,x_n).
\end{equation}
The corresponding eigenvalues now read
\begin{equation}\label{EW}
E^W_\lambda = \sum_{1\leq j\leq n}
\lambda_j\left( \lambda^{}_j +\nu_0+\nu_1+\nu_2+\nu_3 -1+
2(n-j)\nu  \right) .
\end{equation}
\begin{prp}[\cite{die:multivariable}]\label{diffWprp}
The multivariable Wilson polynomials
$p_\lambda^{W}$, $\lambda\in\Lambda$ \eqref{cone},
satisfy the second order difference equation
\begin{equation}\label{eveqW}
D^{W} p_\lambda^{W} \; =\; E_{\lambda}^{W}\: p_\lambda^{W}.
\end{equation}
\end{prp}
For $n=1$, Equation \eqref{eveqW}
reduces to the second order
difference equation for the one-variable Wilson
polynomials \cite{koe-swa:askey-scheme}
\begin{eqnarray}
\frac{ \prod_{0\leq r\leq 3} (i\nu_r+x)}
     { 2ix\, (2ix -1)}
\left( p^{W}_l(x+i)-p^{W}_l(x) \right)&+& \\
\frac{ \prod_{0\leq r\leq 3} (i\nu_r-x)}
     { 2ix\, (2ix +1)}
\left( p^{W}_l(x-i)-p^{W}_l(x) \right) & & \nonumber \\
= l \left( l^{} +\nu_0+\nu_1+\nu_2+\nu_3 -1\right) p^W_l(x) . & & \nonumber
\end{eqnarray}

\subsection{Continuous Hahn type}
For the continuous Hahn type one has
\begin{equation}\label{DcH}
D^{cH} = \sum_{1\leq j\leq n} \left(
V^{cH}_{j,-}(x) (T_{j}-1) \; +\; V^{cH}_{j,+}(x) (T_{j}^{-1}-1) \right)
\end{equation}
with
\begin{eqnarray}
V^{cH}_{j,+}(x) &=& (\nu^+_0 +ix_j) (\nu^+_1 + ix_j)
\prod_{1\leq k \leq n,\, k\neq j}
\left( 1+ \frac{\nu}{i(x_j-x_k)}\right) ,\\
V^{cH}_{j,-}(x) &=& (\nu^-_0 -ix_j) (\nu^-_1 - ix_j)
\prod_{1\leq k \leq n,\, k\neq j}
\left( 1- \frac{\nu}{i(x_j-x_k)}\right) .
\end{eqnarray}
The action of $T_j$ is the same as in the Wilson case (cf. \eqref{Tj})
and the eigenvalues are given by
\begin{equation}\label{EcH}
E^{cH}_\lambda = \sum_{1\leq j\leq n}
\lambda_j\left( \lambda^{}_j +\nu^+_0+\nu^+_1+\nu^-_0+
\nu^-_1 -1+ 2(n-j)\nu  \right) .
\end{equation}
\begin{prp}[\cite{die:multivariable}]\label{diffcHprp}
The multivariable continuous Hahn polynomials
$p_\lambda^{cH}$, $\lambda\in\Lambda$ \eqref{cone},
satisfy the second order difference equation
\begin{equation}\label{eveqcH}
D^{cH} p_\lambda^{cH} \; =\; E_{\lambda}^{cH}\: p_\lambda^{cH}.
\end{equation}
\end{prp}
For $n=1$, Equation \eqref{eveqcH}
reduces to the second order
difference equation for the one-variable continuous Hahn
polynomials \cite{koe-swa:askey-scheme}
\begin{eqnarray}
 (\nu^-_0 -ix) (\nu^-_1 - ix)
\left( p^{cH}_l(x+i)-p^{W}_l(x) \right) &+& \\
(\nu^+_0 +ix) (\nu^+_1 + ix)
\left( p^{cH}_l(x-i)-p^{cH}_l(x) \right)\nonumber  \\
= l \left( l^{} +\nu^+_0+\nu^+_1+\nu^-_0+\nu^-_1 -1\right)
p^{cH}_l(x) . & & \nonumber
\end{eqnarray}

\subsection{Jacobi type}
In the case of multivariable Jacobi type polynomials the operator $D$
diagonalized by $p_\lambda$ is given by a second
order differential operator of the form
\begin{eqnarray}\label{DJ}
\makebox[1em]{} D^J \!\!\!&=&\!\!\! -\sum_{1\leq j\leq n} \partial_j^2 \;\;\;\;
-\alpha \sum_{1\leq j\leq n} \left(
\nu_0 \cot (\frac{\alpha x_j}{2}) -\nu_1 \tan (\frac{\alpha x_j}{2})
\right) \partial_j \\
\!\! & & \!\!\!\!\!-\alpha\nu \!\!\!\!\sum_{1\leq j< k\leq n}
\left( \cot (\frac{\alpha}{2} (x_j+x_k)) (\partial_j +\partial_k) +
\cot (\frac{\alpha}{2} (x_j-x_k)) (\partial_j -\partial_k) \right) \nonumber
\end{eqnarray}
where $\partial_j\equiv \partial /\partial x_j$.
The eigenvalue of $D^J$ on $p^J_\lambda$ takes the value
\begin{equation}
E^J_\lambda = \sum_{1\leq j\leq n}
\lambda_j\left( \lambda^{}_j +\nu_0+\nu_1+
2(n-j)\nu  \right) .
\end{equation}
\begin{prp}[\cite{vre:formulas,deb:systeme}]
The multivariable Jacobi polynomials
$p_\lambda^{J}$, $\lambda\in\Lambda$ \eqref{cone},
satisfy the second order differential equation
\begin{equation}\label{eveqJ}
D^{J} p_\lambda^{J} \; =\; E_{\lambda}^{J}\: p_\lambda^{J}.
\end{equation}
\end{prp}
For $n=1$, Equation \eqref{eveqJ}
reduces to the second order
differential equation for the one-variable Jacobi
polynomials \cite{abr-ste:handbook,koe-swa:askey-scheme}
\begin{eqnarray}
-\frac{d^2p_\lambda^{J}}{dx^2}(x) &-& \alpha \left(
\nu_0 \cot (\frac{\alpha x}{2}) -\nu_1 \tan (\frac{\alpha x}{2})
\right)
\frac{dp_\lambda^{J}}{dx}(x) \\
&=& l \left( l^{} +\nu_0+\nu_1 \right) p^{J}_l(x) . \nonumber
\end{eqnarray}

\section{Limit transitions}\label{sec4}
The operator $D$ of the previous section can be used to arrive at the
following useful representation for the polynomials $p_\lambda$
(cf. \cite{mac:orthogonal,die:self-dual,sto-koo:limit})
\begin{equation}\label{urep}
p_\lambda = \left(
\prod_{\mu \in \Lambda , \mu < \lambda }
\frac{ D-E_\mu}{E_\lambda - E_\mu}
\right) m_\lambda .
\end{equation}
Indeed, it is not difficult to infer that the r.h.s. of \eqref{urep}
determines a polynomial satisfying the defining properties {\em i.} and
{\em ii.} stated in Section~\ref{sec2}. To this end one uses the
Triangularity \eqref{tria} and Symmetry \eqref{sym} of $D$, together
with the observation that in each of the concrete cases discussed
above the denominators in \eqref{urep} are nonzero since
for parameter values indicated in Section~\ref{sec2} one has that
(see \cite[Sec. 5.2]{die:commuting} and \cite{sto-koo:limit})
\begin{equation}
\mu < \lambda \implies E_\mu < E_\lambda .
\end{equation}
(It is immediate from the triangularity of $D$ that
the r.h.s. of \eqref{urep} can be written as a linear combination of
monomials $m_\mu$ with $\mu \leq \lambda$; that the r.h.s.
is also orthogonal to all $m_\mu$ with $\mu < \lambda$
follows from the symmetry of $D$ and the fact that the operator in
the numerator---viz.
$\prod_{\mu \in \Lambda , \mu < \lambda } (D-E_\mu )$---annihilates
the subspace
$\text{Span} \{ m_\mu \}_{\mu\in\Lambda ,\mu <\lambda}$ in view of the
Cayley-Hamilton theorem.)

Below we will use Formula~\eqref{urep} to derive limit transitions
from the Askey-Wilson type to the Wilson, continuous Hahn and Jacobi
type families, respectively. The transition
`Askey-Wilson $\rightarrow$ Jacobi' has already been considered before
in \cite{mac:orthogonal,die:commuting,sto-koo:limit} and is included here
mainly for the sake of completeness. It will be put to use in
Section~\ref{sec6} when deriving a system of recurrence
relations for the multivariable Jacobi type polynomials.

\subsection{Askey-Wilson$\rightarrow$Wilson}
When studying the limit $p_\lambda^{AW}\rightarrow p_\lambda^W$ it is
convenient to first express the multivariable Askey-Wilson polynomials in
terms of a slightly modified monomial basis consisting of
the functions
\begin{equation}
\tilde{m}_\lambda^{AW}(x)= (2/\alpha)^{2|\lambda |}\:
m_{sym,\lambda} (\sin^2(\alpha x_1/2),\ldots ,\sin^2(\alpha x_n/2)),
\;\;\;\; \lambda \in \Lambda ,
\end{equation}
where $|\lambda |\equiv \lambda_1+\cdots +\lambda_n$.
Notice that
\begin{equation}\label{monlimAW-W}
\lim_{\alpha \rightarrow 0}\tilde{m}_\lambda^{AW}(x) = m^W_\lambda (x)
\end{equation}
whereas the original monomials $m^{AW}_\lambda(x)$ \eqref{AWmonom} all reduce
to constant functions in this limit.
Using the relation
$\sin^2 (\alpha x_j/2)=1/2 -(e^{i\alpha x_j}+e^{-i\alpha x_j})/4$
one easily infers that the bases
$\{ \tilde{m}_\lambda^{AW} \}_{\lambda \in\Lambda}$
and $\{ m_\lambda^{AW} \}_{\lambda \in\Lambda}$ are related
by a triangular transformation of the form
\begin{equation}\label{btrafo}
\tilde{m}_\lambda = (-1/\alpha^2)^{|\lambda |}\: m^{AW}_\lambda +
\sum_{\mu \in\Lambda ,\, \mu < \lambda}
a_{\lambda ,\mu}\: m_\mu
\;\;\;\;\;\; \text{with}\;\; a_{\lambda ,\mu}\in {\Bbb R}.
\end{equation}
It is clear that in Formula~\eqref{urep} we may always replace the
monomial basis $\{ m_\lambda\}_{\lambda\in\Lambda}$ by a different basis
that is related by a unitriangular transformation, since
(cf. above) the operator
$\prod_{\mu\in\Lambda ,\mu < \lambda}(D-E_\mu)$ in the numerator
of the r.h.s. annihilates
the subspace
$\text{Span} \{ m_\mu \}_{\mu\in\Lambda ,\mu <\lambda}$ because of the
Cayley-Hamilton theorem. Hence, by taking in account the diagonal matrix
elements in the basis transformation~\eqref{btrafo}, one sees
that Formula~\eqref{urep} can rewritten in terms of
$\tilde{m}_\lambda^{AW}$ as
\begin{equation}\label{urepAW}
p^{AW}_\lambda = (-\alpha^2)^{|\lambda |}\; \left(
\prod_{\mu \in \Lambda , \mu < \lambda }
\frac{ D^{AW}-E^{AW}_\mu}{E^{AW}_\lambda - E^{AW}_\mu}
\right) \tilde{m}^{AW}_\lambda .
\end{equation}
If we now substitute
\begin{equation}\label{AW-Wsub}
q=e^{-\alpha},\;\;\;\; t=e^{-\alpha \nu },\;\;\;\;
t_r=e^{-\alpha \nu_r}\;\; (r=0,1,2,3)
\end{equation}
in $D^{AW}$ \eqref{DAW} and $E_\lambda^W$ \eqref{EAW}, then we have that
\begin{equation}\label{opelimAW-W}
\lim_{\alpha \rightarrow 0} \alpha^{-2}\: D^{AW}  =D^W,
\;\;\;\;\;\;\;\;\;
\lim_{\alpha \rightarrow 0} \alpha^{-2}\: E_\lambda^{AW}
=E_\lambda^W.
\end{equation}
(Notice to this end that for $q=e^{-\alpha}$ the action of
$T_{j,q}$ \eqref{Tjq} on trigonometric polynomials
is the same as that of $T_j$ \eqref{Tj}, i.e., the action amounts
to a shift of the variable $x_j$ over an imaginary unit:
$x_j\rightarrow x_j+i$. To infer then that in the limit
$\alpha \rightarrow 0$ the difference operator $\alpha^{-2}D^{AW}$
formally goes to $D^W$ boils down to checking that the
coefficients of the operator converge as advertised.)

By applying the limits \eqref{monlimAW-W}
and \eqref{opelimAW-W} to
Formula~\eqref{urepAW} we end up with the following limiting
relation between the multivariable Askwey-Wilson and Wilson type
polynomials.
\begin{prp}\label{AW-Wprp}
For Askey-Wilson parameters given by \eqref{AW-Wsub} one has
\begin{equation}
p_\lambda^W (x)=\lim_{\alpha \rightarrow 0}\;
(-1/\alpha^2)^{|\lambda |}\: p_\lambda^{AW}(x),\;\;\;\;\;\;
\lambda\in\Lambda
\end{equation}
(with $|\lambda |\equiv \lambda_1+\cdots +\lambda_n$).
\end{prp}

\subsection{Askey-Wilson$\rightarrow$continuous Hahn}
Just like in the previous subsection, the derivation of the
transition $p_\lambda^{AW}\rightarrow p_\lambda^{cH}$
hinges again on Formula \eqref{urep}.
If we shift the variables $x_1,\ldots ,x_n$ over a half period
by setting
\begin{equation}\label{shift}
x_j \rightarrow x_j - \pi /(2\alpha ),\;\;\;\;\; j=1,\ldots ,n
\end{equation}
and substitute parameters in the following way
\begin{equation}\label{AW-cHsub}
\begin{array}{llll}
q=e^{-\alpha},& t=e^{-\alpha \nu} ,& & \\ [1ex]
t_0=-ie^{-\alpha \nu^+_0},& t_1=-ie^{-\alpha \nu^+_1},&
t_2=ie^{-\alpha \nu^-_0},& t_3=ie^{-\alpha \nu^-_1},
\end{array}
\end{equation}
then the version of Formula \eqref{urep} for the multivariable
Askey-Wilson polynomials takes the form
($e_j$ denotes the $j$th unit vector in the standard basis of
${\Bbb R}^n$)
\begin{equation}\label{urepAW-cH}
p_\lambda^{AW}\left( x-\frac{\pi}{2\alpha}(e_1+\cdots +e_n)\right)
= \prod_{\mu\in\Lambda ,\: \mu <\lambda}
\left( \frac{\tilde{D}^{AW}-E_\mu^{AW}}{E_\lambda^{AW}-E_\mu^{AW}}\right)
\tilde{m}^{AW}_\lambda (x)
\end{equation}
where
\begin{equation*}
\tilde{D}^{AW} = \sum_{1\leq j\leq n} \left(
\tilde{V}^{AW}_j(x) (T_{j}-1) \; +\; \tilde{V}^{AW}_{-j}(x)
(T_{j}^{-1}-1) \right) ,
\end{equation*}
with
\begin{eqnarray*}
\tilde{V}^{AW}_{j}(x) &=&
\frac{ (1+e^{-\alpha \nu^+_0} e^{i\alpha x_j})
       (1+e^{-\alpha \nu^+_1} e^{i\alpha x_j})
       (1-e^{-\alpha \nu^-_0} e^{i\alpha x_j})
       (1-e^{-\alpha \nu^-_1} e^{i\alpha x_j})}
     { (1+e^{2i\alpha x_j})\, (1+e^{-\alpha} e^{2i\alpha x_j})} \\
& & \times \prod_{1\leq k \leq n,\, k\neq j}
\left( \frac{1+e^{-\alpha\nu} e^{i\alpha (x_j +x_k)}}
            {1+e^{i\alpha (x_j+x_k)}}\right)
\left( \frac{1-e^{-\alpha\nu} e^{i\alpha (x_j -x_k)}}
            {1-e^{i\alpha (x_j-x_k)}} \right) \\
\tilde{V}^{AW}_{-j}(x) &=&
\frac{ (1-e^{-\alpha \nu^+_0} e^{-i\alpha x_j})
       (1-e^{-\alpha \nu^+_1} e^{-i\alpha x_j})
       (1+e^{-\alpha \nu^-_0} e^{-i\alpha x_j})
       (1+e^{-\alpha \nu^-_1} e^{-i\alpha x_j})}
     { (1+e^{-2i\alpha x_j})\, (1+e^{-\alpha} e^{-2i\alpha x_j})} \\
& & \times \prod_{1\leq k \leq n,\, k\neq j}
\left( \frac{1+e^{-\alpha\nu} e^{-i\alpha (x_j +x_k)}}
            {1+e^{-i\alpha (x_j+x_k)}}\right)
\left( \frac{1-e^{-\alpha\nu} e^{-i\alpha (x_j -x_k)}}
            {1-e^{-i\alpha (x_j-x_k)}} \right)
\end{eqnarray*}
and
\begin{equation*}
\tilde{m}_\lambda^{AW}(x)\equiv
m_{sym ,\lambda}(2\sin(\alpha x_1),\ldots ,2\sin (\alpha x_n)).
\end{equation*}
(Just as in the case of the
transition Askey-Wilson$\rightarrow$Wilson we have
rewritten the operators $T_{j,q}$ \eqref{Tjq} for $q=e^{-\alpha}$
as $T_j$ \eqref{Tj}.)
After dividing by $(2\alpha)^{|\lambda |}$ the r.h.s. of \eqref{urepAW-cH}
goes for $\alpha\rightarrow 0$ to the corresponding formula for
the continuous Hahn
polynomials (i.e. with $\tilde{D}^{AW}\rightarrow D^{cH}$,
$E_\lambda^{AW}\rightarrow E_\lambda^{cH}$ and
$(2\alpha)^{-|\lambda |} \tilde{m}_\lambda^{AW}\rightarrow m_\lambda^{cH}$).
Hence, we now arrive at the following limiting relation between
the multivariable Askey-Wilson and continuous Hahn type polynomials.
\begin{prp}\label{AW-cHprp}
For Askey-Wilson parameters given by \eqref{AW-cHsub} one has
\begin{equation}
p_\lambda^{cH} (x)=\lim_{\alpha \rightarrow 0}\;
\frac{1}{(2\alpha )^{|\lambda |}}\:
p_\lambda^{AW}\left( x-\frac{\pi}{2\alpha}\omega \right),\;\;\;\;\;\;
\lambda\in\Lambda
\end{equation}
where $\omega \equiv e_1 +\cdots +e_n$ (with $e_j$ denoting the $j$th unit
vector in the standard basis of ${\Bbb R}^n$).
\end{prp}

\subsection{Askey-Wilson$\rightarrow$Jacobi}
To recover the Jacobi type polynomials we substitute the
Askey-Wilson parameters
\begin{equation}\label{AW-Jsub}
t=q^g,\;\;\;\;
t_0=q^{g_0},\;\;\;\;
t_1=-q^{g_1},\;\;\;\;
t_2=q^{g_0^\prime +1/2},\;\;\;\;
t_3=-q^{g_1^\prime +1/2}.
\end{equation}
With these parameters the formula of the Form \eqref{urep}
for the Askey-Wilson type polynomials reduces in the
limit $q\rightarrow 1$
to the corresponding formula for the Jacobi type polynomials
(i.e., the difference operator $D^{AW}$ with eigenvalues
$E^{AW}_\lambda$
gets replaced by the differential operator $D^J$ with eigenvalues
$E^J_\lambda$). The limit $q\rightarrow 1$ amounts to sending the
difference step size to zero. In order to analyze the behavior of
the operator $D^{AW}$ for $q\rightarrow 1$ in detail it is convenient to
substitute $q=e^{-\alpha\beta}$ (so the action of $T_{j,q}$
\eqref{Tjq} on trigonometric polynomials amounts to the shift
$x_j\rightarrow x_j+i\beta$) and then write formally $T_{j,q}=\exp
(i\beta\partial_j)$. A formal expansion in $\beta$ then
shows that $D^{AW}\sim \beta^2 D^J$ and that $E_\lambda^{AW}\sim
\beta^2 E_\lambda^J$ for $\beta \rightarrow 0$. Here the Jacobi
parameters $\nu$, $\nu_r$ are related to the
parameters $g$, $g_r^{(\prime )}$ in \eqref{AW-Jsub} via
$\nu = g$, $\nu_0=g_0+g_0^\prime$ and
$\nu_1=g_1+g_1^\prime$. As a consequence we
obtain the following limiting
relation between the multivariable Askey-Wilson and Jacobi type
polynomials.
\begin{prp}\label{AW-Jprp}
For Askey-Wilson parameters given by \eqref{AW-Jsub} one has
\begin{equation}
p_\lambda^J(x) = \lim_{q\rightarrow 1}\; p_\lambda^{AW}(x)
\end{equation}
with the Jacobi parameters $\nu$, $\nu_0$ and $\nu_1$ taking the value
$g$, $g_0+g_0^\prime$ and $g_1+g_1^\prime$, respectively.
\end{prp}

{\em Remarks: i.} In the above derivations of Propositions \ref{AW-Wprp},
\ref{AW-cHprp} and \ref{AW-Jprp} we have used that the
Askey-Wilson type difference operator converges formally (i.e.,
without specifying the domains of the operators of interest) to the
corresponding operators connected with
the  Wilson, continuous Hahn and Jacobi type polynomials,
respectively. In our case such formal limits get their precise meaning when
being applied to Formula~\eqref{urep}.

{\em ii.} For all our four families AW, W, cH and J the dependence
of (the coefficients of) the operator $D$
and of the eigenvalues $E_\lambda$ is polynomial in the
parameters $t,t_0,t_1,t_2,t_3$ (AW), $\nu,\nu_0,\nu_1,\nu_2,\nu_3$ (W),
$\nu,\nu^\pm_0,\nu^\pm_0$ (cH) and $\nu,\nu_0,\nu_1$ (J), respectively.
Hence it is clear from Formula~\eqref{urep} that (the coefficients of)
the polynomials $p_\lambda$ are rational in these parameters.
We may thus extend the parameter domains for the polynomials
given in Section~\ref{sec2} to generic (complex) values
by alternatively characterizing $p_\lambda$ as
the polynomial of the form
$p_\lambda = m_\lambda + \sum_{\mu\in\lambda ,\mu < \lambda} c_{\lambda ,\mu}
m_\mu$ satisfying the eigenvalue equation
$Dp_\lambda = E_\lambda p_\lambda$.
It is clear that the limit transitions discussed in this
section then extend to these larger parameter domains of
generic (complex) parameter values.

{\em iii.} In the case of one variable the limit transitions from Askey-Wilson
polynomials to Wilson, continuous Hahn and Jacobi polynomials
were collected in \cite{koe-swa:askey-scheme} (together with
many other limits between the various (basic) hypergeometric
orthogonal families appearing in the ($q$-)Askey scheme).

\section{Higher order difference or differential equations}\label{sec5}
In \cite{die:commuting} it was shown that the second order difference
equation for the multivariable Askey-Wilson type polynomials can be
extended to a system of difference equations having the structure of
eigenvalue equations of the form
\begin{equation}\label{diffsystem}
D_r\, p_\lambda = E_{r,\lambda} \, p_\lambda, \;\;\;\;\;
r=1,\ldots ,n,
\end{equation}
for $n$ independent commuting difference operators
$D_1,D_2,\ldots ,D_n$ of order $2,4,\ldots\!,2n$, respectively.
For $r=1$ one recovers the second order difference equation discussed
in Section~\ref{sec3AW}. After recalling the explicit expressions
for the Askey-Wilson type difference
operators $D^{AW}_r$ and their eigenvalues $E^{AW}_{r,\lambda}$,
we will apply the limit transitions of Section~\ref{sec4} to arrive at
similar systems of difference equations
for the multivariable Wilson and continuous Hahn type polynomials.
In case of the transition `Askey-Wilson$\rightarrow$ Jacobi' the
step size is sent to zero and the
system of difference equations degenerates to a system of
hypergeometric differential equations, thus generalizing the state of affairs
for the second order operator in the previous section.
This limit from Askey-Wilson type difference equations
to Jacobi type differential equations has already been discussed in
detail in \cite[Sec. 4]{die:commuting}, so here we will merely state
the results and refrain from presenting a
complete treatment of this case.

\subsection{Askey-Wilson type}
The difference operators diagonalized by
the multivariable Askey-Wilson polynomials via Eq.~\eqref{diffsystem}
are given by
\begin{equation}\label{diffAWr}
D_r^{AW}=\sum_{\stackrel{J\subset \{ 1,\ldots ,n\} ,\, 0\leq |J|\leq r}
               {\varepsilon_j=\pm 1,\; j\in J}}
U^{AW}_{J^c,\, r-|J|}(x)\,  V^{AW}_{\varepsilon J,\, J^c}(x)\,
T_{\varepsilon J,q}
\;\;\;\; r=1,\ldots , n,
\end{equation}
with
\begin{eqnarray*}
T_{\varepsilon J,q} \!\!\! &=&\!\!\! \prod_{j\in J} T_{j,q}^{\varepsilon_j} \\
V^{AW}_{\varepsilon J,\, K}(x)\!\!\! &=&\!\!\!
\prod_{j\in J} w^{AW}(\varepsilon_jx_j)
\prod_{\stackrel{j,j^\prime \in J}{j<j^\prime}}
v^{AW}(\varepsilon_jx_j+\varepsilon_{j^\prime}x_{j^\prime})
v^{AW}(\varepsilon_jx_j+\varepsilon_{j^\prime}x_{j^\prime}
         -i\ln(q)/\alpha )\\
& & \times
\prod_{\stackrel{j\in J}{k\in K}} v^{AW}(\varepsilon_j x_j+x_k)
v^{AW}(\varepsilon_j x_j -x_k),\\
U^{AW}_{K,p}(x)\!\!\! &=&\\ \!\!\!\!\!\!
 (-1)^p\!\!\!\!\!\!   & &\!\!\!\!\!\!\!\!\!
\sum_{\stackrel{L\subset K,\, |L|=p}
               {\varepsilon_l =\pm 1,\; l\in L }}\!\!
\Bigl( \prod_{l\in L} w^{AW}(\varepsilon_l x_l)
\prod_{\stackrel{l,l^\prime \in L}{l<l^\prime}}
v^{AW}(\varepsilon_lx_l+\varepsilon_{l^\prime}x_{l^\prime})
v^{AW}(-\varepsilon_lx_l-\varepsilon_{l^\prime}x_{l^\prime}
+i\ln(q)/\alpha )\\
& &\times
\prod_{\stackrel{l\in L}{k\in K\setminus L}} v^{AW}(\varepsilon_l x_l+x_k)
v^{AW}(\varepsilon_l x_l -x_k) \Bigr) ,
\end{eqnarray*}
and
\begin{eqnarray}\label{vAW}
v^{AW}(z) &=& t^{-1/2} \left( \frac{1-t\, e^{i\alpha z}}
                 {1-e^{i\alpha z}}\right) \\
w^{AW}(z) &=& (t_0 t_1 t_2 t_3 q^{-1})^{-1/2}
\frac{ \prod_{0\leq r\leq 3} (1-t_r\, e^{i\alpha z})}
       {(1-e^{2i\alpha z})
        (1-q\, e^{2i\alpha z})  } .\label{wAW}
\end{eqnarray}
Here the action of the operators
$T_{j,q}^{\pm 1}$ is defined in accordance with
\eqref{Tjq}. The
summation in \eqref{diffAWr} is over all index sets
$J\subset \{ 1,\ldots ,n\}$ with cardinality $|J|\leq r$ and over all
configurations of signs $\varepsilon_j \in \{ +1 ,-1\}$ with $j\in J$.
Furthermore, by convention empty products are taken to be equal to one
and $U_{K,p}\equiv 1$ for $p=0$.

The corresponding eigenvalue of $D^{AW}_r$ on $p^{AW}_\lambda$
has the value
\begin{equation}
E^{AW}_{r,\lambda}=
E_r(\tau_1q^{\lambda_1}+\tau_1^{-1}q^{-\lambda_1},\ldots ,
\tau_nq^{\lambda_n}+\tau_n^{-1}q^{-\lambda_n};
\tau_r +\tau_r^{-1},\ldots ,\tau_n+\tau_n^{-1})
\end{equation}
where
\begin{eqnarray}\label{Er}
& & E_r(x_1,\ldots ,x_n;y_r,\ldots ,y_n) \equiv \\
& & \hspace{7em}
\sum_{\stackrel{ J\subset \{ 1,\ldots ,n\} }{0\leq |J|\leq r}}
(-1)^{r-|J|} \prod_{j\in J} x_j
\sum_{r\leq l_1\leq \cdots \leq l_{r-|J|}\leq n}
y_{l_1}\cdots y_{l_{r-|J|}}  \nonumber
\end{eqnarray}
and
\begin{equation}\label{tauAW}
\tau_j = t^{n-j} (t_0 t_1 t_2 t_3 q^{-1})^{1/2},
\;\;\;\;\;\;\;\;\;\;\;\;\;\;\;\; j=1,\ldots ,n.
\end{equation}
(The second sum in \eqref{Er} is understood to be equal to $1$ when
$|J|=r$.)

Summarizing,
we have the following theorem generalizing Proposition~\ref{diffAWprp}.
\begin{thm}[\cite{die:commuting,die:diagonalization}]\label{diffAWrthm}
The multivariable Askey-Wilson polynomials $p_\lambda^{AW}$,
$\lambda\in\Lambda$ \eqref{cone}
satisfy a system of difference equations of the form
\begin{equation}\label{eveqAWr}
D^{AW}_r\, p_\lambda^{AW}=E_{r,\lambda}^{AW}\, p_\lambda^{AW},
\;\;\;\;\;\; r=1,\ldots ,n.
\end{equation}
\end{thm}
For $r=1$ the Difference equation \eqref{eveqAWr} goes over in
the second order Difference equation \eqref{eveqAW} after multiplication
by a constant with value $t^{n-1}(t_0 t_1 t_2 t_3 q^{-1})^{1/2}$.
More generally, one may multiply the Difference equation \eqref{eveqAWr}
for arbitrary $r$ by the constant factor
$t^{r(n-1)-r(r-1)/2}(t_0 t_1 t_2 t_3 q^{-1})^{r/2}$ to obtain a
difference equation that is polynomial in the parameters
$t$, $t_0,\ldots ,t_3$ and rational in $q$. Such multiplication amounts to
omitting the factors $t^{-1/2}$ and $(t_0 t_1 t_2 t_3 q^{-1})^{-1/2}$
in the definition of
$v^{AW}(z)$ \eqref{vAW} and $w^{AW}(x)$ \eqref{wAW} and to
replacing the eigenvalues by
\begin{eqnarray*}
E^{AW}_{r,\lambda} &\rightarrow & t^{-r(r-1)/2} \\
& & \times E_r(\tau^+_1q^{\lambda_1}+\tau_1^- q^{-\lambda_1},\ldots ,
\tau_n^+q^{\lambda_n}+\tau_n^- q^{-\lambda_n};
\tau_r^+ +\tau_r^- ,\ldots ,\tau_n^++\tau_n^-)
\end{eqnarray*}
with
\begin{eqnarray*}
\tau_j^+&=&t^{n-1}(t_0 t_1 t_2 t_3 q^{-1})^{1/2}\, \tau_j =\;
t_0 t_1 t_2 t_3 q^{-1} t^{2n-1-j}, \\
\tau_j^-&=&t^{n-1}(t_0 t_1 t_2 t_3 q^{-1})^{1/2}\, \tau_j^{-1} =\;
t^{j-1}.
\end{eqnarray*}

\subsection{Wilson type}
If we substitute Askey-Wilson parameters of the form \eqref{AW-Wsub}
and divide by $\alpha^{2r}$, then for $\alpha\rightarrow 0$ the
operator $D^{AW}_r$ goes over in
\begin{equation}\label{diffWr}
D^W_r= \sum_{\stackrel{J\subset \{ 1,\ldots ,n\} ,\, 0\leq |J|\leq r}
               {\varepsilon_j=\pm 1,\; j\in J}}
U^{W}_{J^c,\, r-|J|}(x)\,  V^{W}_{\varepsilon J,\, J^c}(x)\,
T_{\varepsilon J}
\;\;\;\; r=1,\ldots , n,
\end{equation}
with
\begin{eqnarray*}
T_{\varepsilon J}\!\!\! &=&\!\!\! \prod_{j\in J} T_{j}^{\varepsilon_j} \\
V^{W}_{\varepsilon J,\, K}(x)\!\!\! &=&\!\!\!
\prod_{j\in J} w^{W}(\varepsilon_jx_j)
\prod_{\stackrel{j,j^\prime \in J}{j<j^\prime}}
v^{W}(\varepsilon_jx_j+\varepsilon_{j^\prime}x_{j^\prime})
v^{W}(\varepsilon_jx_j+\varepsilon_{j^\prime}x_{j^\prime}+i )\\
& & \times
\prod_{\stackrel{j\in J}{k\in K}} v^{W}(\varepsilon_j x_j+x_k)
v^{W}(\varepsilon_j x_j -x_k),\\
U^{W}_{K,p}(x)\!\!\! &=&\!\!\! \\
\!\!\! (-1)^p\!\!\!\!\!\!\!\!\!    & &\!\!\!\!\!
\sum_{\stackrel{L\subset K,\, |L|=p}
               {\varepsilon_l =\pm 1,\; l\in L }}\;
\Bigl( \prod_{l\in L} w^{W}(\varepsilon_l x_l)
\prod_{\stackrel{l,l^\prime \in L}{l<l^\prime}}
v^{W}(\varepsilon_lx_l+\varepsilon_{l^\prime}x_{l^\prime})
v^{W}(-\varepsilon_lx_l-\varepsilon_{l^\prime}x_{l^\prime}-i )\\
& &\times
\prod_{\stackrel{l\in L}{k\in K\setminus L}} v^{W}(\varepsilon_l x_l+x_k)
v^{W}(\varepsilon_l x_l -x_k) \Bigr) ,
\end{eqnarray*}
and
\begin{equation}\label{vwW}
v^{W}(z) = \frac{ i\nu +z}{z},\;\;\;\;\;\;\;\;
w^{W}(z) =
\frac{ \prod_{0\leq r\leq 3}(i\nu_r+z)}
       {2iz(2iz-1)} .
\end{equation}
Here the action of $T_j^\pm$ is taken to be in accordance with \eqref{Tj}.
To verify this limit it suffices to recall that for $q=e^{-\alpha}$
the action of $T_{j,q}$ \eqref{Tjq} amounts to that of $T_{j}$
\eqref{Tj} and to observe that for parameters \eqref{AW-Wsub}
one has $\lim_{\alpha \rightarrow} v^{AW}(z)=v^W(z)$
and  $\lim_{\alpha \rightarrow} \alpha ^{-2}w^{AW}(z)=w^W(z)$.

The eigenvalues become in this limit
\begin{equation}
E_{r,\lambda}^W=
E_r\left( (\rho_1^W+\lambda_1)^2,\ldots ,(\rho_n^W+\lambda_n)^2;
(\rho^W_r)^2,\ldots ,(\rho^W_n)^2 \right)
\end{equation}
with $E_r(\cdots ;\cdots)$ taken from \eqref{Er} and
\begin{equation}\label{rhoW}
\rho_j^W = (n-j)\nu + (\nu_0+\nu_1+\nu_2+\nu_3-1)/2 .
\end{equation}
For the eigenvalues the computation verifying the limit
is a bit more complicated than for the
difference operators; it hinges on the following lemma.
\begin{lem}[\mbox{\cite[Sec. 4.2]{die:commuting}}]\label{evlim}
One has
\begin{eqnarray*}
&& \lim_{\alpha\rightarrow 0}\;
\alpha^{-2r}
E_r(e^{\alpha x_1}+e^{-\alpha x_1},\ldots ,e^{\alpha x_n}+e^{-\alpha x_n};
    e^{\alpha y_r}+e^{-\alpha y_r},\ldots ,e^{\alpha y_n}+e^{-\alpha y_n}) \\
&&\;\;\;\;\;\;\;\;\;\;\;\;\;\;\;\;\;\;
=E_r(x_1^2,\ldots ,x_n^2;y_r^2,\ldots,y_n^2) . \nonumber
\end{eqnarray*}
\end{lem}
We may thus conclude that the transition AW$\rightarrow$W gives rise to
the following generalization of Proposition~\ref{diffWprp}.
\begin{thm}
The multivariable Wilson polynomials $p_\lambda^{W}$,
$\lambda\in\Lambda$ \eqref{cone}
satisfy a system of difference equations of the form
\begin{equation}\label{eveqWr}
D^{W}_r\, p_\lambda^{W}=E_{r,\lambda}^{W}\, p_\lambda^{W},
\;\;\;\;\;\; r=1,\ldots ,n.
\end{equation}
\end{thm}
For $r=1$ the Difference equation \eqref{eveqWr}
coincides with the second order Difference equation \eqref{eveqW}.

\subsection{Continuous Hahn type}
After shifting over a half period as in \eqref{shift} and choosing
Askey-Wilson parameters of the form \eqref{AW-cHsub}, the operators
$\alpha^{-2r}D^{AW}_r$ go for $\alpha\rightarrow 0$ over in
\begin{equation}\label{diffcHr}
D^{cH}_r= \!\!\!\sum_{\stackrel{J_+,J_-\subset \{ 1,\ldots ,n\}}
               {J_+\cap J_- =\emptyset ,\; |J_+|+|J_-|\leq r}}
\!\!\!\!\!\! U^{cH}_{J_+^c\cap J_-^c,\, r-|J_+|-|J_-|}(x)\,
V^{cH}_{J_+,J_-;\, J_+^c\cap J_-^c}(x)\,
T_{J_+,\, J_-}
\end{equation}
$r=1,\ldots , n$, where
\begin{eqnarray*}
T_{J_+,\,J_-}\!\!\! &=&\!\!\! \prod_{j\in J_+} T_{j}^{-1}
\prod_{j\in J_-} T_{j}\\
V^{cH}_{J_+, J_- ;\, K}(x)\!\!\!  &=&\!\!\!
\prod_{j\in J_+} w^{cH}_+(x_j)
\prod_{j\in J_-} w^{cH}_-(x_j)
\prod_{j\in J_+,j^\prime\in J_-}
v^{cH}(x_j-x_{j^\prime})
v^{cH}(x_j-x_{j^\prime}-i )\\
& & \times
\prod_{\stackrel{j\in J_+}{k\in K}}v^{cH}(x_j-x_k)
\prod_{\stackrel{j\in J_-}{k\in K}}v^{cH}(x_k-x_j),\\
U^{cH}_{K,p}(x)\!\!\! &=&\!\!\! \\
\!\!\!\!\!\!  (-1)^p \!\!\!\!\!\!\!\!\! & &\!\!\!\!\! \!\!\!
\sum_{\stackrel{L_+,L_-\subset K,\, L_+\cap L_-=\emptyset}
               {|L_+|+|L_-|=p}}\;
\Bigl(
\prod_{l\in L_+}\!\! w_+^{cH}(x_l)\!\!
\prod_{l\in L_-}\!\! w_-^{cH}(x_l)\!\!
\prod_{\stackrel{l\in L_+}{l^\prime\in L_-}}\!\!\!\!
v^{cH}(x_l-x_{l^\prime})
v^{cH}(x_{l^\prime}-x_l+i )\\
& &\times
\prod_{\stackrel{l\in L_+}{k\in K\setminus L_+\cup L_-}} \!\!\!\!
v^{cH}(x_l-x_k)\!\!\!\!
\prod_{\stackrel{l\in L_-}{k\in K\setminus L_+\cup L_-}}\!\!\!\!
v^{cH}(x_k -x_l) \Bigr) ,
\end{eqnarray*}
and
\begin{eqnarray}
v^{cH}(z)\!\!\! &=&\!\!\!  (1+\frac{\nu}{iz}) ,\\
w_+^{cH}(z)\!\!\! &=&\!\!\!  (\nu_0^+ +iz)(\nu_1^+ +iz)   ,\;\;\;\;\;\;\;\;\;
w_-^{cH}(z)= (\nu_0^- -iz)(\nu_1^- -iz). \label{wcH}
\end{eqnarray}
(The action of the operators $T_j^\pm$ is again in accordance with \eqref{Tj}.)
To verify this transition one uses that
$v^{AW}(\pm (x_j+x_k))\rightarrow 1$,
$v^{AW}(x_j-x_k)\rightarrow v^{cH}(x_k-x_j)$ and that
$\alpha^{-2}w^{AW}(\pm x_j)\rightarrow w^{cH}_\mp (x_j)$,
if one sends $\alpha$ to zero after having substituted
\eqref{shift} and \eqref{AW-cHsub}.

The computation of the limit of the eigenvalues is exactly the same
as in the Wilson case and the result reads
\begin{equation}
E_{r,\lambda}^{cH}=
E_r\left( (\rho_1^{cH}+\lambda_1)^2,\ldots ,(\rho_n^{cH}+\lambda_n)^2;
(\rho^{cH}_r)^2,\ldots ,(\rho^{cH}_n)^2 \right)
\end{equation}
with $E_r(\cdots ;\cdots)$ taken from \eqref{Er} and
\begin{equation}\label{rhocH}
\rho_j^{cH} = (n-j)\nu + (\nu^+_0+\nu^+_1+\nu^-_0+\nu^-_1-1)/2 .
\end{equation}
Hence, we arrive the following generalization of
Proposition~\ref{diffcHprp}.
\begin{thm}
The multivariable continuous Hahn polynomials $p_\lambda^{cH}$,
$\lambda\in\Lambda$ \eqref{cone}
satisfy a system of difference equations of the form
\begin{equation}\label{eveqcHr}
D^{cH}_r\, p_\lambda^{cH}=E_{r,\lambda}^{cH}\, p_\lambda^{cH},
\;\;\;\;\;\; r=1,\ldots ,n.
\end{equation}
\end{thm}
For $r=1$ the Difference equation \eqref{eveqcHr}
coincides with the second order Difference equation \eqref{eveqcH}.

\subsection{Jacobi type}
Ater substituting Askey-Wilson parameters given by \eqref{AW-Jsub}
and dividing by a constant with value $(1-q)^{2r}$ the $r$th
difference equation
in Theorem~\ref{diffAWrthm} for the Askey-Wilson type polynomials goes
for $q\rightarrow 1$ over in a differential equation
$D^J_r p^J_\lambda = E^J_{r,\lambda}p^J_\lambda$ of order $2r$
for the multivariable Jacobi type polynomials
\cite[Sec. 4]{die:commuting}.
The computation of the eigenvalue
$E^J_{r,\lambda}=\lim_{q\rightarrow 1} (1-q)^{-2r} E_{r,\lambda}^{AW}$
hinges again on Lemma~\ref{evlim} and the result is
\begin{equation}\label{evrJ}
E^J_{r,\lambda} =
E_r\left( (\rho_1^{J}+\lambda_1)^2,\ldots ,(\rho_n^{J}+\lambda_n)^2;
(\rho^{J}_r)^2,\ldots ,(\rho^{J}_n)^2 \right)
\end{equation}
with $E_r(\cdots ;\cdots)$ taken from \eqref{Er} and
\begin{equation}\label{rhoJ}
\rho_j^{J} = (n-j)\nu + (\nu_0+\nu_1)/2 ,
\end{equation}
where $\nu =g$, $\nu_0=g_0+g_0^\prime$ and $\nu_1=g_1+g_1^\prime$.
For $r=1$ the differential operator
$D^J_r=\lim_{q\rightarrow 1} (1-q)^{-2r} D_r^{AW}$ is given by
$D^J$ \eqref{DJ}. More generally one has that $D^J_r$ is of the form
\begin{equation}\label{DJr}
D^J_r = \sum_{\stackrel{J\subset \{ 1,\ldots ,n\}}{|J|=r}}
\prod_{j\in J} \partial_j^2 \;\;\; + \text{l.o.}
\end{equation}
(where l.o. stands for the parts of lower order in the partials),
but it seems difficult to obtain the relevant differential
operators for arbitrary $r$ in explicit form starting from
$D^{AW}_r$ \eqref{diffAWr}.
\begin{thm}[\mbox{\cite[Sec. 4]{die:commuting}}]
The multivariable Jacobi polynomials $p_\lambda^{J}$,
$\lambda\in\Lambda$ \eqref{cone}
satisfy a system of differential equations of the form
\begin{equation}\label{eveqJr}
D^{J}_r\, p_\lambda^{J}=E_{r,\lambda}^{J}\, p_\lambda^{J},
\;\;\;\;\;\; r=1,\ldots ,n,
\end{equation}
where $D^J_r=\lim_{q\rightarrow 1} (1-q)^{-2r} D_r^{AW}$ is of the
form \eqref{DJr} and the corresponding eigenvalues are given by
\eqref{evrJ}.
\end{thm}
For $r=1$ the Differential equation \eqref{eveqJr}
coincides with the second order Differential equation \eqref{eveqJ}.

{\em Remarks: i.} The proof in
\cite{die:commuting,die:diagonalization} demonstrating
that the multivariable Askey-Wilson type polynomials
satisfy the system of difference equations in
Theorem~\ref{diffAWrthm} runs along the same lines as the proof
for the special case when $r=1$ (Proposition~\ref{diffAWprp}):
it consists of demonstrating
that the operator $D^{AW}_r$ is triangular with respect to the
monomial basis $\{ m^{AW}_\lambda \}_{\lambda\in \Lambda}$
and that it is symmetric with respect to the inner product
$\langle \cdot ,\cdot \rangle_{\Delta^{AW}}$.
The triangularity proof consists of two parts. First it is shown
that $D^{AW}_r m_\lambda^{AW}$ lies in the space ${\cal H}^{AW}$
by inferring that poles originating from the zeros in the
denominators of $v^{AW}(z)$ \eqref{vAW} and $w^{AW}(z)$ \eqref{wAW}
all cancel each other.
Next, it is verified that the operator $D_r^{AW}$ indeed
subtriangular by analyzing the asymptotics for $x$ at infinity.

It is noteworthy to observe that if one tries to apply the same
approach to arrive at a direct proof of the system of difference
equations for the multivariable Wilson and continuous Hahn type
polynomials
(i.e., without using the fact that these two cases may be seen as limiting
cases of the Askey-Wilson type polynomials), then some complications
arise. In both cases
the proof that the operator $D_r$ maps
the space ${\cal H}$ into itself and that
it is symmetric with respect to
the inner product $\langle \cdot ,\cdot \rangle_{\Delta}$
applies without significant changes. However, it is now not so
simple to deduce from the asymptotics at infinity
that the difference operator $D_r$ is indeed triangular.
The reason---in a nutshell---is that
the functions $w^W(z)$ \eqref{vwW} and $w_\pm^{cH}(z)$ \eqref{wcH}
no longer have constant asymptotics for $z\rightarrow \infty$
(as does $w^{AW}(z)$ for $\exp(i\alpha z) \rightarrow \infty$).
As a consequence, it is now much more difficult than
for the Askey-Wilson case to rule out the a priori possibility
that monomials $m_\mu$ with $\mu\not\leq \lambda$ enter the
expansion of $D_r m_\lambda$ in monomial symmetric functions
(in principle $D_r$ might a priori raise the degrees of the polynomials).
For $r=1$ one easily checks by inspection that such $m_\mu$
with $\mu \not\leq \lambda$ indeed do not appear in the
expansion of $D_r m_\lambda$ but for general $r$ this is not
so easily seen from the explicit expressions at hand.

{\em ii.}
Some years ago, a rather explicit characterization of a family of
commuting differential operators simultaneously diagonalized by
the Jacobi type polynomials
$p_\lambda^J$ (and generating the same algebra as the operators
$D^J_1,\ldots ,D_n^J$)
was presented by Debiard \cite{deb:systeme}.
Alternatively, it also turned out possible to
express such differential operators in terms of
symmetric functions of Heckman's trigonometric
generalization of the Dunkl differential-reflection operators
related to the root system $BC_n$ \cite{hec:elementary,dun:differential}.
In both cases, however, it is a nontrivial problem
to deduce from those results an explicit combinatorial formula for
(the coefficients of) the operator $D_r^J$ for arbitrary $r$.
In the case of Debiard's operators one problem is that
the corresponding
eigenvalues do not seem to be known precisely in closed form
(and also that one would like to commute all coefficients to the left);
in the case of an expression in terms of trigonometric
Dunkl type differential-reflection operators
it appears to be difficult to explicitly compute the differential operator
corresponding to the restriction of the relevant symmetrized
differential-reflection operators to the space
of symmetric polynomials (except when the order of the symbol is small).
In the latter case the problem is that one has to commute all
reflection operators to the right (and preferably
all coefficients to the left).
This poses a combinatorial exercise that
seems tractable only for small order of the symbol.

{\em iii.} If one transforms the operator $D^J=D^J_1$ to a second order
differential operator that is self-adjoint in
$L^2({\Bbb R}^n,d x_1\cdots d x_n)$ by conjugating with the square
root of the weight function $\Delta^J(x)$ \eqref{weightJ},
then one arrives at a Schr\"odinger operator
of the form
\begin{eqnarray}\makebox[2em]{}
(\Delta^J)^{-\frac{1}{2}} D^J (\Delta^J)^{\frac{1}{2}} = &&\\
&& \hspace{-10em}-\sum_{1\leq j\leq n} \partial_j^2
\;\;+\;\; \frac{1}{4}\: \alpha^2\sum_{1\leq j\leq n}\Bigl(
\frac{\nu_0 (\nu_0-1)}{\sin^2 (\frac{\alpha x_j}{2})} +
\frac{\nu_1 (\nu_1-1)}{\cos^2 (\frac{\alpha x_j}{2})}\Bigr) \nonumber \\
&& \hspace{-8em} +\frac{1}{2}\: \nu (\nu-1)\: \alpha^2
\sum_{1\leq j< k\leq n} \Bigl(
\frac{1}{\sin^2 \frac{\alpha}{2}(x_j-x_k)} +
\frac{1}{\sin^2 \frac{\alpha}{2}(x_j-x_k)}\Bigr)  \;\;\;\; -\varepsilon_0.
\nonumber
\end{eqnarray}
where $\varepsilon_0= \alpha^2\sum_{1\leq j\leq n} (\rho_j^J)^2$.
Explicit expressions for $n$ independent commuting differential operators
generating the same commutative algebra as transformed operators
$(\Delta^J)^{-\frac{1}{2}} D^J_r (\Delta^J)^{\frac{1}{2}}$,
$r=1,\ldots ,n$ can be found in the literature
as a special case of the formulas
presented in \cite{och-osh-sek:commuting,osh-sek:commuting}.
However, just as in the previous remark
it is again nontrivial to determine explicitly the
relation between our transformed operators
$(\Delta^J)^{-\frac{1}{2}} D^J_r (\Delta^J)^{\frac{1}{2}}$,
$r=1,\ldots ,n$ and the relevant specialization of
the differential operators in
\cite{och-osh-sek:commuting,osh-sek:commuting},
because to our knowledge the eigenvalues of the latter operators
are not available in closed form (except in cases when the
order of the symbol is small).

\section{Recurrence relations}\label{sec6}
In this section we first recall the
system of recurrence
relations for the multivariable
Askey-Wilson type polynomials presented in
\cite{die:self-dual}. Next, the limit transitions of
Section~\ref{sec4} are applied to arrive at similar systems of
recurrence relations for the multivariable Wilson, continuous Hahn and
Jacobi type polynomials, respectively.
To describe these recurrence relations it is convenient to pass from
the monic polynomials $p_\lambda(x)$
to a different normalization by introducing
\begin{equation}\label{reno}
P_\lambda (x) \equiv c_\lambda\: p_\lambda (x),
\;\;\;\;\;\;\;\;\;\;\;\;\;\;
c_\lambda = c^{|\lambda |}\:
\frac{\hat{\Delta}_+(\rho )}{ \hat{\Delta}_+(\rho +\lambda )} ,
\end{equation}
where $c$ denotes some constant not depending on $\lambda$ (recall
also that $|\lambda |\equiv \lambda_1 +\cdots +\lambda_n$) and
the function $\hat{\Delta}_+ (x)$ is of the form
\begin{equation}\label{deltahp}
 \hat{\Delta}_+ (x) =\prod_{1\leq j<k \leq n }
\hat{d}_{v,+} (x_j+x_{k})\, \hat{d}_{v,+} (x_j-x_{k})\;
\prod_{1\leq j\leq n}
\hat{d}_{w,+} (x_j).
\end{equation}
The precise value of the constant $c$, the vector
$\rho =(\rho_1,\ldots ,\rho_n)$ and the form of the
functions $\hat{d}_{v,+}$, $\hat{d}_{w,+}$ depends on the type of polynomials
of interest and will be detailed below separately for each case.

The general structure of the recurrence relations for the renormalized
polynomials $P_\lambda (x)$ reads
\begin{equation}\label{recr}
\hat{E}_r (x)\: P_{\lambda} (x)=
\sum_{\stackrel{J\subset \{ 1,\ldots ,n\} ,\, 0\leq|J|\leq r}
               {\varepsilon_j=\pm 1,\; j\in J;\;
                e_{\varepsilon J} +\lambda \in \Lambda}}
\!\!\!\!\!\!\!\!\!
\hat{U}_{J^c,\, r-|J|}(\rho +\lambda)\,
\hat{V}_{\varepsilon J,\, J^c}(\rho +\lambda)\,
P_{\lambda +e_{\varepsilon J}} (x)
\end{equation}
with $r=1,\ldots ,n$, and
\begin{eqnarray*}
e_{\varepsilon J} &= & \sum_{j\in J} \varepsilon_j e_j ,\\
\hat{V}_{\varepsilon J,\, K}(x)\!\!\! &=&\!\!\!
\prod_{j\in J} \hat{w}(\varepsilon_jx_j)
\prod_{\stackrel{j,j^\prime \in J}{j<j^\prime}}
\hat{v}(\varepsilon_jx_j+\varepsilon_{j^\prime}x_{j^\prime})\,
\hat{v}(\varepsilon_jx_j+\varepsilon_{j^\prime}x_{j^\prime}
         +1)\\
& & \times
\prod_{\stackrel{j\in J}{k\in K}} \hat{v}(\varepsilon_j x_j+x_k)\,
\hat{v}(\varepsilon_j x_j -x_k),\\
\hat{U}_{K,p}(x)\!\!\! &=&\\ \!\!\!\!\!\!
 (-1)^p\!\!\!\!\!\!   & &\!\!\!\!\!\!\!\!\!
\sum_{\stackrel{L\subset K,\, |L|=p}
               {\varepsilon_l =\pm 1,\; l\in L }}\!\!
\Bigl( \prod_{l\in L} \hat{w}(\varepsilon_l x_l)
\prod_{\stackrel{l,l^\prime \in L}{l<l^\prime}}
\hat{v}(\varepsilon_lx_l+\varepsilon_{l^\prime}x_{l^\prime})\,
\hat{v}(-\varepsilon_lx_l-\varepsilon_{l^\prime}x_{l^\prime}
-1 )\\
& &\times
\prod_{\stackrel{l\in L}{k\in K\setminus L}}
\hat{v}(\varepsilon_l x_l+x_k)\,
\hat{v}(\varepsilon_l x_l -x_k) \Bigr) .
\end{eqnarray*}
The functions $\hat{E}_1(x),\ldots ,\hat{E}_n(x)$ appearing in the
l.h.s. of \eqref{recr} denote certain (explicitly given)
symmetric polynomials
that generate the algebra of all symmetric polynomials.
The $r$th recurrence relation expresses the fact that
the expansion of the product
$\hat{E}_r(x) P_{\lambda} (x)$
in terms of the basis elements
$P_{\mu} (x)$, $\mu\in\Lambda$ is known explicitly (i.e., the
coefficients in the expansion are known in closed form).
The expansion coefficients are determined by the functions $\hat{v}$
and $\hat{w}$ (together with the vector $\rho$)
whose precise form again
depends on the class of polynomials
of interest.
In combinatorics one sometimes refers to this type of recurrence relations
as to generalized Pieri type formulas after similar expansion formulas for the
products of elementary or complete symmetric functions and
Schur functions (in terms of the latter functions) \cite{mac:symmetric}.
In the simplest case, i.e. when $r=1$, the Recurrence formula
\eqref{recr} reduces to an expression of the form
\begin{eqnarray}\label{rec1}
\hat{E}(x) P_\lambda (x) &=&
\sum_{\stackrel{1\leq j\leq n}{\lambda + e_j\in\Lambda}}
\hat{V}_j (\rho +\lambda)
\left( P_{\lambda+e_j}(x)-P_\lambda (x) \right) +\\
& &\sum_{\stackrel{1\leq j\leq n}{\lambda - e_j\in\Lambda}}
\hat{V}_{-j} (\rho +\lambda)
\left( P_{\lambda-e_j}(x)-P_\lambda (x)\right)\nonumber
\end{eqnarray}
with $\hat{E}(x)=\hat{E}_1(x)$ and
\begin{equation}\label{vjh}
\hat{V}_{\pm j}(x) =
\hat{w}(\pm x_j) \prod_{1\leq k\leq n,\, k\neq j}
\hat{v}(\pm x_j +x_k)\, \hat{v}(\pm x_j -x_k).
\end{equation}
We will see that for the four families considered below
this formula reduces in the case of one variable to the well-known
three-term recurrence relation for the Askey-Wilson, Wilson,
continuous Hahn or Jacobi polynomials, respectively.

{\em Note:}
To date a complete  proof (contained in \cite{die:self-dual}) of the
recurrence relations for the multivariable Askey-Wilson type
polynomials exists only for parameters satisfying the
(self-duality) condition
\begin{equation}\label{self-dualc}
q\, t_0\, t_1^{-1} t_2^{-1} t_3^{-1} =1 .
\end{equation}
It is for this reason that in the present formulation
of the system of recurrence relations for our multivariable
hypergeometric polynomials given below
some restrictions on the parameters are imposed (except in the case
of Jacobi type polynomials where such a restriction turns out to
be not necessary).
As is argued in \cite{die:self-dual} for the Askey-Wilson
type, however, it is expected
(and easily checked for the special case of one variable) that
our recurrence relations remain valid also for parameter values
not meeting these restrictions (cf. Remark {\em i.}
at the end of this section).

\subsection{Askey-Wilson type}
The normalization constants $c^{AW}_\lambda$ for the
polynomials in \eqref{reno} of Askey-Wilson type
$P^{AW}_\lambda(x)$ are determined by the constant
$c^{AW}=1$ and the functions
\begin{eqnarray}
\hat{d}^{AW}_{v,+}(z)& =&t^{-z/2}
\frac{(q^{z}; q)_\infty}
     {( t q^{z}; q)_\infty},  \label{dhAWvp}  \\[1ex]
\hat{d}^{AW}_{w,+}(z)&=& (\hat{t}_0 \hat{t}_1 \hat{t}_2 \hat{t}_3q^{-1})^{-z/2}
\frac{(q^{2 z};\,
                                    q)_\infty    }
     {(\hat{t}_0 q^{z},\,
       \hat{t}_1 q^{z},\,
       \hat{t}_2 q^{z},\,
       \hat{t}_3 q^{z};\,
                                    q)_\infty    }  .\label{dhAWwp}
\end{eqnarray}
Here we have introduced dependent
parameters $\hat{t}_r$ that are related to the
Askey-Wilson parameters $t_r$ by
\begin{equation}\label{hattr}
\begin{array}{lll}
\hat{t}_0 &=& (t_0 t_1      t_2      t_3      q^{-1})^{1/2},\\
\hat{t}_1 &=& (t_0 t_1      t_2^{-1} t_3^{-1} q     )^{1/2},\\
\hat{t}_2 &=& (t_0 t_1^{-1} t_2      t_3^{-1} q     )^{1/2},\\
\hat{t}_3 &=& (t_0 t_1^{-1} t_2^{-1} t_3      q     )^{1/2}.
\end{array}
\end{equation}
The vector $\rho =\rho^{AW}=(\rho^{AW}_1,\ldots ,\rho_n^{AW})$ has
in the present case the components
\begin{equation}\label{rhoAW}
\rho_j^{AW} = {}^q\log \tau_j ,\;\;\;\;\;\;\;
\tau_j =t^{n-j} (t_0 t_1 t_2 t_3 q^{-1})^{1/2}
\end{equation}
and is introduced mostly for notational
convenience. In fact, the logarithm with base $q$
entering through the components
$\rho^{AW}_j$ merely has a formal meaning and appears in our
formulas always as an exponent of $q$.

The symmetric functions $\hat{E}^{AW}_r (x)$
multiplying $P^{AW}_{\lambda} (x)$ in
the l.h.s. of \eqref{recr} are given by
\begin{equation}\label{ErhAW}
\hat{E}^{AW}_r (x) =
E_r (e^{i\alpha x_1}+e^{-i\alpha x_1},\ldots ,
     e^{i\alpha x_n}+e^{-i\alpha x_n};
\hat{\tau}_r+\hat{\tau}_r^{-1},\ldots ,
\hat{\tau}_n+\hat{\tau}_n^{-1} )
\end{equation}
with $E_r(\cdots ;\cdots)$ being taken from \eqref{Er} (as usual) and with
\begin{equation}\label{tauhAW}
\hat{\tau}_j=t^{n-j}
(\hat{t}_0 \hat{t}_1 \hat{t}_2 \hat{t}_3 q^{-1})^{1/2}
\end{equation}
(cf. \eqref{rhoAW}).
The coefficients in the r.h.s. of \eqref{recr}
entering the expansion of the function
$\hat{E}^{AW}_r(x) P^{AW}_{\lambda}(x)$ in terms of the basis elements
$P^{AW}_{\mu}(x)$ are characterized by the functions
\begin{eqnarray}\label{vhAW}
\hat{v}^{AW}(z) &=& t^{-1/2} \left( \frac{ 1-t\, q^{z}}
                 { 1-q^{z}} \right) ,\\
\hat{w}^{AW}(z) &=&
(\hat{t}_0\hat{t}_1\hat{t}_2\hat{t}_3 q^{-1})^{-1/2}
\frac{ \prod_{0\leq r\leq 3}
(1-\hat{t}_r\, q^{z})}
       {(1-q^{2z})
        (1-q^{2z+1})  } \label{whAW}
\end{eqnarray}
(cf. \eqref{vAW}, \eqref{wAW}).
We now have the following theorem from \cite{die:self-dual}.
\begin{thm}[\cite{die:self-dual}]\label{recrAWthm}
The renormalized
multivariable Askey-Wilson type polynomials $P_\lambda^{AW}(x)$,
$\lambda \in \Lambda$ \eqref{cone}, satisfy a system of
recurrence relations given by \eqref{recr} (with
$\hat{E}_r^{AW}$, $\hat{v}^{AW}$,
$\hat{w}^{AW}$ and $\rho^{AW}$ taken from \eqref{ErhAW},
\eqref{vhAW}, \eqref{whAW} and \eqref{rhoAW}, and with $r=1,\ldots ,n$)
for Askey-Wilson parameters subject
to the condition $q\, t_0\, t_1^{-1} t_2^{-1} t_3^{-1} =1$.
\end{thm}

For $r=1$ one has
$$\hat{E}^{AW}(x)=\hat{E}_1^{AW}(x)= \sum_{1\leq j\leq n}
\left( e^{i\alpha x_j}+ e^{-i\alpha x_j}
-\hat{\tau}_j-\hat{\tau}_j^{-1} \right) . $$
The corresponding Recurrence formula \eqref{rec1}, \eqref{vjh}
coincides in the case of one variable with the
three-term recurrence relation for the renormalized
Askey-Wilson polynomials
\begin{equation}\label{P1AW}
P^{AW}_l(x) =
{}_4\phi_3
\left(
\begin{matrix}
q^{-l},\; t_0 t_1 t_2 t_3 q^{l-1},\;
t_0 e^{i\alpha x},\; t_0 e^{-i\alpha x} \\ [0.5ex]
t_0 t_1 ,\; t_0 t_2 ,\; t_0 t_3
\end{matrix} \; ; q,q \right) ,
\end{equation}
which reads \cite{ask-wil:some,koe-swa:askey-scheme}
\begin{eqnarray}
& & (2\cos (\alpha x) -t_0-t_0^{-1})P^{AW}_l(x) = \\
& &\frac{(1-t_0t_1t_2t_3q^{l-1})\prod_{1\leq r\leq 3}(1-t_0t_rq^l)}
     {t_0(1-t_0t_1t_2t_3q^{2l-1})(1-t_0t_1t_2t_3q^{2l})}
\left( P^{AW}_{l+1}(x)-P^{AW}_l(x) \right) + \nonumber \\
& &\frac{t_0 (1-q^{l})\prod_{1\leq r<s\leq 3}(1-t_st_rq^{l-1})}
     {(1-t_0t_1t_2t_3q^{2l-2})(1-t_0t_1t_2t_3q^{2l-1})}
\left( P^{AW}_{l-1}(x)-P^{AW}_l(x) \right) .\nonumber
\end{eqnarray}
It is clear that in this special situation Theorem~\ref{recrAWthm} indeed holds
without restriction on the parameters as was anticipated more
generally for arbitrary $r$ and $n$ in \cite{die:self-dual}.

It is important to observe that, despite the appearances of square
roots (and infinite products) in intermediate expressions,
both the normalization constants $c^{AW}_\lambda$ and
the Askey-Wilson type Recurrence relations \eqref{recr} are in the end
rational in $q$ and the parameters $t$, $t_0,\ldots ,t_3$.
For the recurrence relations this is is rather immediate from the
fact that
\begin{eqnarray}
&&\hat{v}^{AW}(\varepsilon_j(\rho^{AW}_j+\lambda_j)+\varepsilon_k
(\rho^{AW}_k+\lambda_k))=
t^{-1/2} \left(
\frac{ 1-t\tau_j^{\varepsilon_j}\tau_k^{\varepsilon_k}
q^{\varepsilon_j\lambda_j+ \varepsilon_K\lambda_k}}
{ 1-\tau_j^{\varepsilon_j}\tau_k^{\varepsilon_k}
q^{\varepsilon_j\lambda_j+ \varepsilon_K\lambda_k}} \right)\\
&&\hat{w}^{AW}(\varepsilon_j(\rho^{AW}_j+\lambda_j))=
(\hat{t}_0\hat{t}_1\hat{t}_2\hat{t}_3 q^{-1})^{-1/2}
\frac{ \prod_{0\leq r\leq 3}
(1-\hat{t}_r\, \tau^{\varepsilon_j}q^{\varepsilon_j\lambda_j})}
       {(1-\tau^{2\varepsilon_j}q^{2\varepsilon_j\lambda_j}
        (1-\tau^{2\varepsilon_j}q^{2\varepsilon_j\lambda_j+1})  }.
\end{eqnarray}
The point is that the combinations
$\tau_j^{\varepsilon_j}\tau_k^{\varepsilon_k}$,
$\hat{t}_r\, \tau^{\varepsilon_j}$ and $\tau^{2\varepsilon_j}$
are rational in $t$, $t_0,\ldots ,t_3$ plus that
the functions $\hat{v}^{AW}$ always emerge in pairs
in \eqref{recr} and
$(\hat{t}_0\hat{t}_1\hat{t}_2\hat{t}_3 q^{-1})^{1/2}=t_0$
(hence the square root constant factors in
$\hat{v}^{AW}$, $\hat{w}^{AW}$ do not spoil the rationality).
The last equality is also needed to see that $\hat{E}_r^{AW}(x)$
\eqref{ErhAW} in the l.h.s. of the recurrence relation
depends rationally on the parameters:
$\hat{\tau}_j= t^{n-j}(\hat{t}_0\hat{t}_1\hat{t}_2\hat{t}_3
q^{-1})^{1/2}= t^{n-j}t_0$.
For $c^{AW}_\lambda$ the rationality in the parameters
is seen similarly after rewriting in the form
(by cancelling common factors in numerator and denominator)
\begin{eqnarray}\label{crew}
c_\lambda^{AW} &=& C
\prod_{1\leq j<k\leq n}
\frac{(\tau_j\tau_k ;\,q)_{\lambda_j+\lambda_k}}
     {(t \tau_j\tau_k ;\,q)_{\lambda_j+\lambda_k}}
\frac{(\tau_j\tau_k^{-1} ;\,q)_{\lambda_j-\lambda_k}}
     {(t \tau_j\tau_k^{-1} ;\,q)_{\lambda_j-\lambda_k}} \\
& & \times \prod_{1\leq j\leq n}
\frac{(\tau_j^2;\, q)_{2\lambda_j}}
     {(\hat{t}_0\tau_j,\, \hat{t}_1\tau_j,\, \hat{t}_2\tau_j,\,
       \hat{t}_3\tau_j;\, q)_{\lambda_j}}  \nonumber
\end{eqnarray}
where $C=\prod_{1\leq j\leq n} \hat{\tau}_j^{\lambda_j}$.

\subsection{Wilson type}
If we substitute Askey-Wilson parameters in accordance
with \eqref{AW-Wsub}, then
the polynomials $P^{AW}_\lambda(x)$ converge for
$\alpha\rightarrow 0$ to
renormalized multivariable Wilson type polynomials
$P^W_\lambda(x)=c_\lambda^W p_\lambda (x)$.
The normalization constant $c_\lambda^W$ is of the form in
\eqref{reno}, \eqref{deltahp}
with $c^W=-1$, the vector $\rho^W$ taken from
\eqref{rhoW}, and with the functions $\hat{d}^{W}_{v,+}$,
$\hat{d}^{W}_{w,+}$
given by
\begin{equation}\label{dhWvwp}
\hat{d}^{W}_{v,+}(z)=
\frac{\Gamma (\nu +z)}
     {\Gamma (z)},    \;\;\;\;\;\;\;\;\;\;\;\;\;\;
\hat{d}^{W}_{w,+}(z)=
\frac{\prod_{0\leq r\leq 3}\Gamma (\hat{\nu}_r +z)   }
     {\Gamma (2z) }  .
\end{equation}
Here we have introduced dependent parameters $\hat{\nu}_r$
related to the Wilson parameters $\nu_r$ by
\begin{equation}\label{nuhr}
\begin{array}{lll}
\hat{\nu}_0 &=& (\nu_0+ \nu_1 +\nu_2 +\nu_3 -1)/2,\\
\hat{\nu}_1 &=& (\nu_0+ \nu_1 -\nu_2 -\nu_3 +1)/2,\\
\hat{\nu}_2 &=& (\nu_0- \nu_1 +\nu_2 -\nu_3 +1)/2,\\
\hat{\nu}_3 &=& (\nu_0- \nu_1 -\nu_2 +\nu_3 +1)/2.
\end{array}
\end{equation}
To verify that the Askey-Wilson type polynomials $P^{AW}_\lambda (x)$
with parameters \eqref{AW-Wsub} indeed converge for $\alpha
\rightarrow 0$ to the
renormalized Wilson type polynomials $P^W_\lambda (x)$ thus defined, one
uses Proposition~\ref{AW-Wprp} and the fact that (for parameters
\eqref{AW-Wsub}) $$\lim_{\alpha\rightarrow 0}
(-\alpha^2)^{|\lambda |} c_\lambda^{AW}= c_\lambda^W.$$
The latter limit can be easily checked with the aid
of Representation \eqref{crew} for $c_\lambda^{AW}$ entailing
\begin{eqnarray*}
\lim_{\alpha\rightarrow 0}\; (-\alpha^2)^{|\lambda |}c_\lambda^{AW}
\!\! &=&\!\! (-1)^{|\lambda |}
\prod_{1\leq j<k\leq n}
\frac{(\rho^W_j+\rho^W_k)_{\lambda_j+\lambda_k}}
     {(\nu +\rho^W_j+\rho^W_k)_{\lambda_j+\lambda_k}}
\frac{(\rho^W_j-\rho^W_k)_{\lambda_j-\lambda_k}}
     {(\nu+ \rho^W_j-\rho^W_k)_{\lambda_j-\lambda_k}} \\
& & \times \prod_{1\leq j\leq n}
\frac{(2\rho^W_j)_{2\lambda_j}}
     {(\hat{\nu}_0+\rho^W_j,\, \hat{\nu}_1+\rho^W_j,\, \hat{\nu}_2+\rho^W_j,\,
       \hat{\nu}_3+\rho^W_j)_{\lambda_j}} . \nonumber
\end{eqnarray*}
It is not difficult to see that the r.h.s. of this expression
indeed coincides with the
above defined $c_\lambda^{W}$ of the form \eqref{reno}, \eqref{deltahp}
by rewriting all Pochhammer symbols
as a quotient of two gamma functions.

Let us next turn to the limiting behavior of the recurrence relations.
For parameters as in \eqref{AW-Wsub} and after division by
$\alpha^{2r}$ the $r$th recurrence relation for the Askey-Wilson type
polynomial $P^{AW}_\lambda (x)$ goes over in a recurrence relation
of the form in \eqref{recr}
for the Wilson type polynomial $P^{W}_\lambda (x)$.
The relevant symmetric function in the l.h.s. of the resulting
Wilson type recurrence relation is given by
\begin{equation}\label{ErhW}
\hat{E}^{W}_r (x) =(-1)^{r}
\sum_{\stackrel{ J\subset \{ 1,\ldots ,n\} }{0\leq |J|\leq r}}
 \prod_{j\in J} x_j^2
\sum_{r\leq l_1\leq \cdots \leq l_{r-|J|}\leq n}
(\hat{\rho}^W_{l_1}\cdots \hat{\rho}^W_{l_{r-|J|}})^2
\end{equation}
with
\begin{equation}
\hat{\rho}^W_j = (n-j)\nu +
(\hat{\nu}_0+\hat{\nu}_1+\hat{\nu}_2+\hat{\nu}_3-1)/2
\end{equation}
(cf. \eqref{rhoW}).
Furthermore, the coefficients in the r.h.s. of this recurrence
relation are determined by the
vector $\rho^W$ \eqref{rhoW} and the functions
\begin{equation}\label{vwhW}
\hat{v}^{W}(z) = \frac{ \nu + z}{ z},\;\;\;\;\;\;\;\;\;
\hat{w}^{W}(z) =
\frac{ \prod_{0\leq r\leq 3} (\hat{\nu}_r +z)}{2z(2z+1)}
\end{equation}
(cf. \eqref{vwW}).
To check this transition from the Askey-Wilson to the Wilson type
recurrence relations one uses that
for Askey-Wilson parameters given by \eqref{AW-Wsub} one has
$$\lim_{\alpha\rightarrow 0}\; \alpha^{-2r} \hat{E}_r^{AW}(x)=
\hat{E}_r^W(x)$$ (in view of Lemma~\ref{evlim}) and that
\begin{eqnarray*}
& &\lim_{\alpha\rightarrow 0}\;
\hat{v}^{AW}(\varepsilon_j \rho^{AW}_j+\varepsilon_k\rho^{AW}_k+z)=
\hat{v}^W(\varepsilon_j\rho^{W}_j+\varepsilon_k\rho^{W}_k+ z) ,\\
& &\lim_{\alpha\rightarrow 0}\;
\alpha^{-2}\hat{w}^{AW}(\varepsilon_j\rho_j^{AW}+z)=
\hat{w}^W(\varepsilon_j\rho_j^W+z).
\end{eqnarray*}

We have thus derived the following theorem.
\begin{thm}\label{recrWthm}
The renormalized
multivariable Wilson type polynomials $P_\lambda^{W}(x)$,
$\lambda \in \Lambda$ \eqref{cone}, satisfy a system of
recurrence relations given by \eqref{recr}  (with
$\hat{E}_r^{W}$, $\hat{v}^{W}$,
$\hat{w}^{W}$ and $\rho^W$ taken from \eqref{ErhW},
\eqref{vwhW} and \eqref{rhoW}, and with $r=1,\ldots ,n$)
for Wilson parameters subject
to the condition $\nu_0 -\nu_1 -\nu_2 -\nu_3 +1=0$.
\end{thm}
The condition on the parameters in Theorem~\ref{recrWthm} is an
immediate consequence of the parameter condition in
Theorem~\ref{recrAWthm} and may be omitted once we are sure
that Theorem~\ref{recrAWthm} holds for general parameters.

For $r=1$ the recurrence formula is given by \eqref{rec1}, \eqref{vjh}
with
\begin{equation}
\hat{E}^{W}(x)=\hat{E}_1^W(x)= -\sum_{1\leq j\leq n}
\left( x_j^2+(\hat{\rho}^W_j)^2 \right) .
\end{equation}
In the case of one variable this recurrence formula coincides with the
three-term recurrence relation for the renormalized
Wilson polynomials
\begin{equation}
P^{W}_l(x) =
{}_4F_3
\left(
\begin{matrix}
-l,\; \nu_0+\nu_1+\nu_2+\nu_3+l-1,\;
\nu_0+ix,\; \nu_0-ix \\ [0.5ex]
\nu_0+\nu_1 ,\; \nu_0+\nu_2 ,\; \nu_0+\nu_3
\end{matrix} \; ; 1\right) ,
\end{equation}
which reads \cite{koe-swa:askey-scheme}
\begin{eqnarray}
& & \makebox[1em]{} -( x^2 +\nu_0^2 )P^{W}_l(x) = \\
& &\frac{(l+\nu_0+\nu_1+\nu_2+\nu_3-1)\prod_{1\leq r\leq 3}(l+\nu_0+\nu_r)}
     {(2l+\nu_0+\nu_1+\nu_2+\nu_3-1)(2l+\nu_0+\nu_1+\nu_2+\nu_3)}
\left( P^{W}_{l+1}(x)-P^{W}_l(x) \right) + \nonumber \\
& &\frac{l\; \prod_{1\leq r<s\leq 3}(l+\nu_s+\nu_r-1)}
        {(2l+\nu_0+\nu_1+\nu_2+\nu_3-2)(2l+\nu_0+\nu_1+\nu_2+\nu_3-1)}
\left( P^{W}_{l-1}(x)-P^{W}_l(x) \right) .\nonumber
\end{eqnarray}
This three-term recurrence formula indeed holds without restriction on
the parameters.

\subsection{Continuous Hahn type}
If we substitute Askey-Wilson parameters given by \eqref{AW-cHsub}
and shift the variables $x_j$ over a half period (cf. \eqref{shift}),
then for $\alpha\rightarrow 0$ the Askey-Wilson
type polynomials $P_\lambda^{AW}(x)$ become renormalized multivariable
continuous Hahn
type polynomials $P_\lambda^{cH}(x)=c_\lambda^{cH}p_\lambda^{cH}(x)$
of the form in \eqref{reno}, \eqref{deltahp} with
$c^{cH}=1/i$, the vector $\rho^{cH}$ taken from \eqref{rhocH}, and
with the functions $\hat{d}^{cH}_{v,+}$, $\hat{d}^{cH}_{w,+}$ reading
\begin{equation}\label{dhcHvwp}
\hat{d}^{cH}_{v,+}(z)=
\frac{\Gamma (\nu +z)}
     {\Gamma (z)},    \;\;\;\;\;\;\;\;\;\;\;\;\;\;
\hat{d}^{cH}_{w,+}(z)=
\frac{\prod_{0\leq r\leq 2}\Gamma (\hat{\nu}_r +z)   }
     {\Gamma (2z) }  .
\end{equation}
Here we have again employed dependent the parameters
$\hat{\nu}_0,\ldots ,\hat{\nu}_3$ that
are now related to the Wilson parameters $\nu^\pm_0,\nu^\pm_1 $ by
\begin{equation}\label{nuhrcH}
\begin{array}{lll}
\hat{\nu}_0 &=& (\nu^+_0+ \nu^+_1 +\nu^-_0 +\nu^-_1 -1)/2,\\
\hat{\nu}_1 &=& (\nu^+_0- \nu^+_1 +\nu^-_0 -\nu^-_1 +1)/2,\\
\hat{\nu}_2 &=& (\nu^+_0- \nu^+_1 -\nu^-_0 +\nu^-_1 +1)/2,\\
\hat{\nu}_3 &=& (\nu^+_0+ \nu^+_1 -\nu^-_0 -\nu^-_1 +1)/2.
\end{array}
\end{equation}
(Notice, however, that $\hat{d}^{cH}_{w,+}$ depends only on
$\hat{\nu}_0,\hat{\nu}_1 ,\hat{\nu}_2$ and not on $\hat{\nu}_3$,
which is merely introduced for convenience and will be
used below (cf. \eqref{rhohcH}).)
To verify the transition
$P_\lambda^{AW}(x)\rightarrow P_\lambda^{cH}(x)$
one uses Proposition~\ref{AW-cHprp} and
the fact that (for parameters \eqref{AW-cHsub})
$\lim_{\alpha\rightarrow 0}
(2\alpha)^{|\lambda |} c_\lambda^{AW}= c_\lambda^{cH}$.
The derivation of the limit $(2\alpha)^{|\lambda |}
c_\lambda^{AW}\rightarrow c_\lambda^{cH}$ is very similar
to that of the Wilson case and hinges again on Representation
\eqref{crew}: first one writes $c^{AW}_\lambda$ \eqref{crew}
explicitly as a rational expression in the
Askey-Wilson parameters $t$, $t_0,\ldots ,t_3$ by invoking of the
definitions \eqref{hattr} and \eqref{rhoAW}; next, substituting
the parameters \eqref{AW-cHsub} and sending $\alpha$ to zero after having
multiplied by $(2\alpha)^{|\lambda |}$ entails
an expression for $c_\lambda^{cH}$ involving Pochhammer symbols
\begin{eqnarray*}
c_\lambda^{cH}&=&(1/i)^{|\lambda |}
\prod_{1\leq j<k\leq n}
\frac{(\rho^{cH}_j+\rho^{cH}_k)_{\lambda_j+\lambda_k}}
     {(\nu +\rho^{cH}_j+\rho^{cH}_k)_{\lambda_j+\lambda_k}}
\frac{(\rho^{cH}_j-\rho^{cH}_k)_{\lambda_j-\lambda_k}}
     {(\nu+ \rho^{cH}_j-\rho^{cH}_k)_{\lambda_j-\lambda_k}} \\
& & \times \prod_{1\leq j\leq n}
\frac{(2\rho^{cH}_j)_{2\lambda_j}}
{(\hat{\nu}_0+\rho^{cH}_j,\,\hat{\nu}_1+\rho^{cH}_j,\,
 \hat{\nu}_2+\rho^{cH}_j)_{\lambda_j}} , \nonumber
\end{eqnarray*}
which is seen to be equal to the stated expression of the form
\eqref{reno}, \eqref{deltahp} by
rewriting all Pochhammer symbols as a quotient of gamma functions.

The corresponding recurrence relations for
the renormalized multivariable Wilson type
polynomials $P^{cH}_\lambda (x)$
are of the form given by \eqref{recr} with the symmetric functions
in the l.h.s. given by
\begin{equation}\label{ErhcH}
\hat{E}^{cH}_r (x) = (-1)^r
\sum_{\stackrel{ J\subset \{ 1,\ldots ,n\} }{0\leq |J|\leq r}}
 \prod_{j\in J} ix_j
\sum_{r\leq l_1\leq \cdots \leq l_{r-|J|}\leq n}
\hat{\rho}^{cH}_{l_1}\cdots \hat{\rho}^{cH}_{l_{r-|J|}}
\end{equation}
with
\begin{equation}\label{rhohcH}
\hat{\rho}^{cH}_j =(n-j)\nu +
(\hat{\nu}_0+\hat{\nu}_1+\hat{\nu}_2+\hat{\nu}_3-1)/2 .
\end{equation}
The coefficients in the r.h.s. are determined by the vector
$\rho^{cH}$ \eqref{rhocH} and the functions
\begin{equation}\label{vwhcH}
\hat{v}^{cH}(z) = \frac{ \nu + z}{ z},\;\;\;\;\;\;\;\;\;
\hat{w}^{cH}(z) =
\frac{ \prod_{0\leq r\leq 2} (\hat{\nu}_r +z)}{2z(2z+1)} .
\end{equation}
To arrive at these recurrence relations for the multivariable
continuous Hahn type polynomials we have substituted
the parameters \eqref{AW-cHsub},
shifted the variables over a half period \eqref{shift} and
divided both sides by $(2i\alpha)^{r}$. For $\alpha\rightarrow 0$
the $r$th Askey-Wilson type recurrence relation then passes over
in the $r$th continuous Hahn type recurrence relation.
At this point one uses that (for parameters
as in \eqref{AW-cHsub})
$$\lim_{\alpha \rightarrow 0}\;\;  (2\alpha i)^{-r}
\hat{E}_r^{AW}(x-\frac{\pi}{2\alpha}))
=\hat{E}^{cH}_r(x)$$ (which immediate from the definition of
$\hat{E}_r^{AW}(x)$) and that
\begin{eqnarray*}
& &\lim_{\alpha\rightarrow 0}\;
\hat{v}^{AW}(\varepsilon_j \rho^{AW}_j+\varepsilon_k \rho^{AW}_k+z)=
\hat{v}^{cH}(\varepsilon_j\rho^{cH}_j+\varepsilon_k\rho^{cH}_k+ z) ,\\
& &\lim_{\alpha\rightarrow 0}\;  (2i\alpha)^{-1}
\hat{w}^{AW}(\varepsilon_j\rho_j^{AW}+z)=
\hat{w}^{cH}(\varepsilon_j\rho_j^{cH}+z).
\end{eqnarray*}

We arrive at the following theorem.
\begin{thm}\label{recrcHthm}
The renormalized
multivariable continuous Hahn type polynomials $P_\lambda^{cH}(x)$,
$\lambda \in \Lambda$ \eqref{cone}, satisfy a system of
recurrence relations given by \eqref{recr} (with
$\hat{E}_r^{cH}$, $\hat{v}^{cH}$,
$\hat{w}^{cH}$ and $\rho^{cH}$ taken from \eqref{ErhcH},
\eqref{vwhcH} and \eqref{rhocH}, and with $r=1,\ldots ,n$)
for continuous Hahn  parameters subject
to the condition $\nu_0^+ -\nu^+_1 -\nu^-_0 -\nu^-_1 +1=0$.
\end{thm}
The condition on the continuous Hahn
parameters $\nu_0^\pm ,\nu_1^\pm$ in Theorem~\ref{recrcHthm} is of
course an artifact of the parameter condition in
Theorem~\ref{recrAWthm} and may again
be omitted as soon as it is shown that
Theorem~\ref{recrAWthm} holds for general parameters.

For $r=1$ the recurrence relation is
of the form in \eqref{rec1}, \eqref{vjh}
with
\begin{equation}
\hat{E}^{cH}(x)= \hat{E}_1^{cH}(x)=-\sum_{1\leq j\leq n}
\left( ix_j+\hat{\rho}^{cH}_j \right) .
\end{equation}
In the case of one variable this formula coincides with the
three-term recurrence relation for the renormalized
continuous Hahn polynomials
\begin{equation}\label{P1W}
P^{cH}_l(x) =
{}_3F_2
\left(
\begin{matrix}
-l,\; \nu^+_0+\nu^-_0+\nu^+_1+\nu^-_1+l-1,\;
\nu^+_0+ix \\ [0.5ex]
\nu^+_0+\nu^-_0 ,\; \nu^+_0+\nu^-_1
\end{matrix} \; ; 1\right) ,
\end{equation}
which reads \cite{koe-swa:askey-scheme}
\begin{eqnarray}
& & \makebox[2em]{} -( ix+\nu_0^+ )P^{cH}_l(x) = \\
& &\frac{(l+\nu^+_0+\nu^+_1+\nu^-_0+\nu^-_1-1)
         (l+\nu_0^+ +\nu_0^-) (l+\nu_0^+ +\nu_1^-)    }
     {(2l+\nu^+_0+\nu^+_1+\nu^-_0+\nu^-_1-1)
      (2l+\nu^+_0+\nu^+_1+\nu^-_0+\nu^-_1)}
\left( P^{cH}_{l+1}(x)-P^{cH}_l(x) \right) + \nonumber \\
& &\frac{-l (l+\nu_1^+ +\nu_0^- -1) (l+\nu_1^+ +\nu_1^- -1)}
        {(2l+\nu^+_0+\nu^+_1+\nu^-_0+\nu^-_1-2)
         (2l+\nu^+_0+\nu^+_1+\nu^-_0+\nu^-_1-1)}
\left( P^{cH}_{l-1}(x)-P^{cH}_l(x) \right) \nonumber
\end{eqnarray}
and holds for general parameters.

\subsection{Jacobi type}
After substituting the Askey-Wilson parameters \eqref{AW-Jsub}
and sending $q$ to one, the Askey-Wilson
type polynomials $P_\lambda^{AW}(x)$ pass over to renormalized
multivariable Jacobi
type polynomials $P_\lambda^{J}(x)=c^J_\lambda \: p_\lambda^J(x)$
of the form in \eqref{reno}, \eqref{deltahp} with the Jacobi parameters
$\nu_0$ and $\nu_1$ taking the value $g_0+g_0^\prime$ and
$g_1+g_1^\prime$, respectively. The normalization constant
$c_\lambda^J$ is determined by the constant
$c^{J}=2^{-2}$ and the vector $\rho^{J}$ taken from \eqref{rhoJ}
together with
the functions $\hat{d}^{J}_{v,+}$ and $\hat{d}^{J}_{w,+}$ given by
\begin{equation}\label{dhJvwp}
\hat{d}^{J}_{v,+}(z)=
\frac{\Gamma (\nu +z)}
     {\Gamma (z)},    \;\;\;\;\;\;\;\;\;\;\;\;\;\;
\hat{d}^{J}_{w,+}(z)=
\frac{\Gamma (\hat{\nu}_0 +z) \Gamma (\hat{\nu}_1 +z)  }
     {\Gamma (2z) }  .
\end{equation}
Here we have introduced dependent parameters $\hat{\nu}_0$,
$\hat{\nu}_1$ related to the Jacobi parameters $\nu_0$, $\nu_1$ by
\begin{equation}\label{nuhrJ}
\begin{array}{lll}
\hat{\nu}_0 &=& (\nu_0+ \nu_1)/2,\\
\hat{\nu}_1 &=& (\nu_0- \nu_1+1)/2.\\
\end{array}
\end{equation}
The verification of the transition
$P^{AW}_\lambda (x)\rightarrow P^{J}_\lambda (x)$
hinges on Proposition~\ref{AW-Jprp} and
the fact that (for parameters \eqref{AW-Jsub})
$\lim_{q\rightarrow 1}
c_\lambda^{AW}= c_\lambda^{J}$.
The latter limit is again checked using Formula \eqref{crew} for
$c_\lambda^{AW}$ entailing for $q\rightarrow 1$
\begin{eqnarray*}
c_\lambda^{J}&=& 2^{-2|\lambda |}
\prod_{1\leq j<k\leq n}
\frac{(\rho^{J}_j+\rho^{J}_k)_{\lambda_j+\lambda_k}}
     {(\nu +\rho^{J}_j+\rho^{J}_k)_{\lambda_j+\lambda_k}}
\frac{(\rho^{J}_j-\rho^{J}_k)_{\lambda_j-\lambda_k}}
     {(\nu+ \rho^{J}_j-\rho^{J}_k)_{\lambda_j-\lambda_k}} \\
& & \times \prod_{1\leq j\leq n}
\frac{(2\rho^{J}_j)_{2\lambda_j}}
{(\hat{\nu}_0+\rho^{J}_j,\,\hat{\nu}_1+\rho^{J}_j )_{\lambda_j}} , \nonumber
\end{eqnarray*}
which seen to be equal to the stated expression for $c_\lambda^J$ of
the form \eqref{reno}, \eqref{deltahp}
by once again rewriting the Pochhammer symbols in terms of gamma functions.

The corresponding recurrence relations for $P^{J}_\lambda (x)$
are of the form in \eqref{recr}, with the
symmetric function in the l.h.s. given by
\begin{equation}\label{ErhJ}
\hat{E}^{J}_r (x) =
(-1)^r \sum_{\stackrel{ J\subset \{ 1,\ldots ,n\} }{|J|= r}}
 \prod_{j\in J} \sin^2\left( \frac{\alpha x_j}{2}\right)
\end{equation}
(this is $(-1)^r$ times the $r$th elementary symmetric function in the
variables $\sin^2(\frac{\alpha x_j}{2})$, $j=1,\ldots ,n$) and
the coefficients in the r.h.s. being determined by
$\rho^{J}$ \eqref{rhoJ} and the functions
\begin{equation}\label{vwhJ}
\hat{v}^{J}(z) = \frac{ \nu + z}{ z},\;\;\;\;\;\;\;\;\;
\hat{w}^{J}(z) =
\frac{ (\hat{\nu}_0 +z) (\hat{\nu}_1 +z)}{2z(2z+1)} .
\end{equation}
One obtains the recurrence relations for the multivariable
Jacobi type polynomials for $q\rightarrow 1$
from the Askey-Wilson type recurrence relations
with parameters \eqref{AW-Jsub} if one divides
both sides by $2^{2r}$. To check this one uses that (for parameters
given by \eqref{AW-Jsub})
$$\lim_{q \rightarrow 1}\; 2^{-2r}\hat{E}_r^{AW}(x)
=\hat{E}^{J}_r(x)$$
(which is rather immediate from the definition of $\hat{E}_r^{AW}(x)$) and that
\begin{eqnarray*}
& &\lim_{q\rightarrow 1}\;
\hat{v}^{AW}(\varepsilon_j \rho^{AW}_j+ \varepsilon_k\rho^{AW}_k+z)=
\hat{v}^{J}(\varepsilon_j \rho^{J}_j+  \varepsilon_k\rho^{J}_k+ z) ,\\
& &\lim_{q\rightarrow 1}\;
2^{-2}\hat{w}^{AW}(\varepsilon_j\rho_j^{AW}+z)=
\hat{w}^{J}(\varepsilon_j\rho_j^{J}+z).
\end{eqnarray*}

We thus arrive at the following theorem.
\begin{thm}\label{recrJthm}
The renormalized
multivariable Jacobi type polynomials $P_\lambda^{J}(x)$,
$\lambda \in \Lambda$ \eqref{cone}, satisfy a system of
recurrence relations given by \eqref{recr}
(with $\hat{E}_r^{J}$, $\hat{v}^{J}$,
$\hat{w}^{J}$ and $\rho^{J}$ taken from \eqref{ErhJ},
\eqref{vwhJ} and \eqref{rhoJ}, and with $r=1,\ldots ,n$).
\end{thm}
Notice that in the Jacobi case it is not needed to impose any condition on
the parameters $\nu_0$, $\nu_1$. The point is that if we substitute the
Askey-Wilson parameters
\eqref{AW-Jsub} then the condition
$q\, t_0\, t_1^{-1} t_2^{-1} t_3^{-1} =1$ in Theorem~\ref{recrAWthm}
gives rise to the condition $g_0-g_0^\prime - g_1-g_1^\prime =0$
on the parameters $g_0$, $g_1$, $g_0^\prime$, $g_1^\prime$.
However, given this condition the Jacobi parameters
$\nu_0=g_0+g_0^\prime$ and $\nu_1=g_1+g_1^\prime$ can still take
arbitrary values. In other words: the confluence of the parameters in
the limit $q\rightarrow 1$ has as consequence that the full
three-parameter family of multivariable Jacobi type
polynomials  (parametrized by $\nu$, $\nu_0$ and $\nu_1$) may be
seen as limiting case of the family of
(self-dual) Askey-Wilson type polynomials
with parameters satisfying \eqref{self-dualc}.

For $r=1$ we now have a recurrence formula of the form in \eqref{rec1}
with
\begin{equation}
\hat{E}^{J}(x)=\hat{E}_1^{J}(x)= -\sum_{1\leq j\leq n}
\sin^2 \left( \frac{\alpha x_j}{2}\right) .
\end{equation}
In the case of one variable this recurrence formula reduces to the
three-term recurrence relation for the renormalized
Jacobi polynomials
\begin{equation}
P^{J}_l(x) =
{}_2F_1
\left(
\begin{matrix}
-l,\; \nu_0+\nu_1+l \\ [0.5ex]
\nu_0+1/2 \end{matrix} \; ; \; \sin^2 \left( \frac{\alpha x}{2}\right) \right)
\end{equation}
reading \cite{abr-ste:handbook,koe-swa:askey-scheme} (notice,
however, that our
normalization of the polynomials differs slightly from the standard
normalization)
\begin{eqnarray}
& & -\sin^2 \left( \frac{\alpha x}{2}\right) \: P^{J}_l(x) = \\
& &\frac{(l+\nu_0+\nu_1) (l+\nu_0 +1/2)}
     {(2l+\nu_0+\nu_1) (2l+\nu_0+\nu_1+1)}
\left( P^{J}_{l+1}(x)-P^{J}_l(x) \right) + \nonumber \\
& &\frac{l (l+\nu_1 -1/2)}
        {(2l+\nu_0+\nu_1) (2l+\nu_0+\nu_1-1)}
\left( P^{J}_{l-1}(x)-P^{J}_l(x) \right) .\nonumber
\end{eqnarray}

{\em Remarks: i.} The combinatorial structure
of the recurrence relations for the multivariable Askey-Wilson
type polynomials is very similar to that of
the difference equations in Section~\ref{sec5}.
This is by no means a coincidence.
In fact, in \cite{die:self-dual} the recurrence relations were derived
from the difference equations
with the aid of a duality property for the
renormalized multivariable
Askey-Wilson polynomials first conjectured by Macdonald
\begin{equation}\label{dualAW}
P_\lambda^{AW}\left(\frac{i\ln (q)}{\alpha}
(\hat{\rho}^{AW}+\mu)\right)
= \hat{P}_\mu^{AW}\left(\frac{i\ln (q)}{\alpha} (\rho^{AW}+\lambda
)\right) .
\end{equation}
Here $\hat{P}_\mu^{AW}(x)$, $\mu\in\Lambda$, denotes the renormalized
multivariable Askey-Wilson type
polynomial with the parameters $t_0,\ldots ,t_3$ being replaced by
the parameters $\hat{t}_0,\ldots ,\hat{t}_3$ (cf. \eqref{hattr})
and $\rho^{AW}$ is given by \eqref{rhoAW} whereas
$\hat{\rho}^{AW}$ is the corresponding vector with
$\tau_j$ \eqref{rhoAW} replaced by $\hat{\tau}_j$ \eqref{tauhAW}.
If one substitutes $x=\frac{i\ln (q)}{\alpha} (\rho^{AW}+\lambda )$
in the $r$th difference equation of Theorem~\ref{diffAWrthm} for
the polynomial $\hat{P}_\mu^{AW}(x)$, then by using
Property~\eqref{dualAW} (and the fact that the coefficients have a
zero at  $x=\frac{i\ln (q)}{\alpha} (\rho^{AW}+\lambda )$ if
$\lambda + e_{\varepsilon J}\neq \Lambda $ \eqref{cone}) one arrives at the
$r$th recurrence relation of Theorem~\ref{recrAWthm} (first
at the points $x=\frac{i\ln (q)}{\alpha} (\hat{\rho}^{AW}+\mu)$,
$\mu\in\Lambda$ and then for arbitrary $x$ using the fact that one deals
with an equality between trigonometric polynomials).
In \cite{die:self-dual} we proved Property~\eqref{dualAW} (and hence
the recurrence relations) for parameters satisfying
the self-duality condition \eqref{self-dualc}
(implying that $\hat{t}_r=t_r$). It is of course expected, however, that
Macdonald's conjecture \eqref{dualAW}
holds for general parameters. In the case of one variable \eqref{dualAW}
this is immediate from the explicit expression of the polynomials
in terms of the terminating basic hypergeometric series \eqref{P1AW}.

{\em ii.} By applying the transition Askey-Wilson $\rightarrow$ Wilson
to \eqref{dualAW} one
arrives at a similar duality relation for the multivariable
Wilson type variables
\begin{equation}\label{dualW}
P_\lambda^{W}\left(i(\hat{\rho}^{W}+\mu)\right) =
\hat{P}_\mu^{W}\left( i(\rho^{W}+\lambda )\right) ,
\end{equation}
where $\hat{P}_\mu^W(x)$ now denotes the Wilson type polynomial with
the parameters $\nu_0,\ldots ,\nu_3$ being replaced by
$\hat{\nu}_0,\ldots ,\hat{\nu}_3$ (cf. \eqref{nuhr}).
The self-duality condition for the parameters inherited from
\eqref{self-dualc} then becomes $\nu_0-\nu_1-\nu_2-\nu_3 +1=0$
(cf. Theorem~\ref{recrWthm}), which implies that $\hat{\nu}_r=\nu_r$.
For $n=1$ the relation \eqref{dualW} is again immediate for general
parameters from the explicit expression of the polynomials in
terms of the terminating hypergeometric series
\eqref{P1W}.

{\em iii.} For $\mu =0$ the r.h.s. of \eqref{dualAW} becomes identical to
one and one arrives at the relation
$P^{AW}_\lambda (i\alpha^{-1} \ln \hat{\tau} )=1$ (where
$\ln \hat{\tau}$ stands for the vector
$(\ln\hat{\tau}_1,\ldots ,\ln\hat{\tau}_n)$ with $\hat{\tau}_j$
taken from \eqref{tauhAW}).
This equality amounts to the following evaluation or specialization
formula for the monic Askey-Wilson type polynomials
(cf. Definition \eqref{reno})
\begin{equation}\label{specfAW}
p_\lambda^{AW}(i\alpha^{-1} \ln \tau ) =
\frac{\hat{\Delta}^{AW}_+(\rho^{AW} +\lambda)}
     {\hat{\Delta}^{AW}_+(\rho^{AW})}.
\end{equation}
Using the limit transitions of Section~\ref{sec4} one
arrives at similar evaluation formulas for the multivariable
Wilson, continuous Hahn and Jacobi type polynomials:
\begin{eqnarray}
p_\lambda^{W}(i\hat{\rho}^W )& =& (-1)^{|\lambda |}\;
\frac{\hat{\Delta}^{W}_+(\rho^{W} +\lambda)}
     {\hat{\Delta}^{W}_+(\rho^{W})}, \\
p_\lambda^{cH}(i\hat{\rho}^{cH} )& =& i^{|\lambda |}\;
\frac{\hat{\Delta}^{cH}_+(\rho^{cH} +\lambda)}
     {\hat{\Delta}^{cH}_+(\rho^{cH})}, \\
p_\lambda^{J}(0)& =& 2^{2|\lambda |}\;
\frac{\hat{\Delta}^{J}_+(\rho^{J} +\lambda)}
     {\hat{\Delta}^{J}_+(\rho^{J})}.
\end{eqnarray}
For the Askey-Wilson case this evaluation formula was conjectured
by Macdonald and then proven
by Cherednik \cite{che:macdonalds}
for special parameters related to
reduced root systems (and admissible pairs of the form $(R,R^\vee)$)
with the aid of shift operators.
A similar proof for the Jacobi case can be found in \cite{opd:some}.
The proof of the Askey-Wilson type
specialization formula \eqref{specfAW}
based on the recurrence relations for the more general
parameters subject to the condition \eqref{self-dualc} was first presented
in \cite{die:self-dual}.

\section{Orthogonality and normalization}\label{sec7}
The difference/differential equations in Section~\ref{sec5}
hold for generic complex parameter values in view of the
rational dependence on the parameters (cf. also Remark {\em ii.} of
Section~\ref{sec4}). The same is true for
the recurrence relations in Section~\ref{sec6} (although we
still have to
impose the parameter restrictions stated in the theorems of
Section~\ref{sec6} so as to keep our present derivation of
the recurrence relations to remain valid, cf. Remark {\em i} of
Section~\ref{sec6}).
In order to interpret the polynomials as an orthogonal
system with weight function $\Delta$, however,
rather than allowing generic complex parameters
we will from now on
always choose the parameters from the domains given in
Section~\ref{sec2}. (The conditions on the parameters
in Section~\ref{sec2} ensure that the relevant
weight functions $\Delta (x)$ are
positive and that the integrals defining the associated
inner products
$\langle \cdot ,\cdot \rangle_\Delta$ converge in absolute value.)

It is immediate from the definition in Section~\ref{sec2}
that the polynomial $p_\lambda$ is orthogonal to $p_\mu$ for $\mu < \lambda$.
It turns out that also for weights that are not comparable with
respect to the partial order \eqref{po} the associated polynomials
are orthogonal.

\begin{thm}[Orthogonality]\label{orthothm}
The multivariable Askey-Wilson, Wilson, continuous Hahn and Jacobi
type polynomials $p_\lambda$, $\lambda\in\Lambda$, form
an orthogonal system with respect to the $L^2$ inner product
with weight function $\Delta (x)$, i.e.
\begin{equation}
\langle p_\lambda ,p_\mu \rangle_\Delta=0
\;\;\;\;\;\;\text{if}\;\;\;\;\;\; \lambda \neq \mu
\end{equation}
(for parameters with values in the domains given in Section~\ref{sec2}).
\end{thm}
The proof of this theorem hinges on the following proposition.
\begin{prp}[Symmetry]\label{symmetry}
The operators $D_1,\ldots,D_n$ of Section~\ref{sec5} are symmetric
with respect to the inner product $\langle \cdot , \cdot
\rangle_\Delta$, i.e.
\begin{equation}
\langle D_r m_\lambda , m_\mu \rangle_\Delta =
\langle  m_\lambda , D_r m_\mu \rangle_\Delta
\end{equation}
(for parameters with values in the domains given in Section~\ref{sec2}).
\end{prp}
A detailed proof of the symmetry for the Askey-Wilson type difference
operators can be found in \cite[Sec. 3.4]{die:commuting}.
The same reasoning used there
can also be applied to prove the proposition for the
Wilson and continuous Hahn case. The Jacobi case follows from the
Askey-Wilson case with the aid of the limit transition from Askey-Wilson type
to Jacobi type polynomials and operators (we refer to
\cite[Sec. 4]{die:commuting} for the precise details).

By combining Proposition~\ref{symmetry}
with the results of Section~\ref{sec5} we conclude
that the polynomials $p_\lambda$, $\lambda\in \Lambda$, are
joint eigenfunctions of $n$ independent operators
$D_1,\ldots,D_n$ that are symmetric with respect to
$\langle \cdot ,\cdot \rangle_\Delta$. The corresponding (real) eigenvalues
$E_{1,\lambda}, \ldots ,E_{n,\lambda}$
separate the points of the integral cone $\Lambda$ \eqref{cone}, i.e.,
if $E_{r,\lambda}=E_{r,\mu}$ for $r=1,\ldots ,n$ then
$\lambda$ must be equal to $\mu$ (this is seen using the fact that the
functions $E_r(x_1,\ldots ,x_n;y_r,\ldots y_n)$ \eqref{Er} generate
the algebra of permutation symmetric polynomials in
the variables $x_1,\ldots ,x_n$).
It thus follows that the polynomials $p_\lambda$ and $p_\mu$ with
$\lambda \neq \mu$ are orthogonal
with respect to the inner product
$\langle \cdot ,\cdot \rangle_\Delta$ as eigenfunctions
of a symmetric operator corresponding to different eigenvalues.
Another proof for the orthogonality of the multivariable
Askey-Wilson type polynomials was given by Koornwinder in \cite{koo:askey}.
In the case of Jacobi type polynomials independent orthogonality
proofs (not using the fact that the Jacobi type is a degenerate
case of the Askey-Wilson type) can be found
in e.g. \cite{deb:systeme} or in \cite{hec:elementary} (upon
specialization to the root system $BC_n$).

The main purpose of this section is to express the squared norms
of the polynomials (viz. $\langle p_\lambda ,p_\lambda \rangle_\Delta$)
in terms of the squared norm of the unit polynomial
(viz. $\langle 1 , 1 \rangle_\Delta$, which corresponds to $\lambda =0$).
Our main tool to achieve this goal
consists of the recurrence relations derived in the
preceding section.
Starting point is the identity
\begin{equation}\label{ident}
\langle \hat{E}_r P_\lambda , P_{\lambda +\omega_r}\rangle_\Delta
=
\langle  P_\lambda , \hat{E}_r P_{\lambda +\omega_r}\rangle_\Delta ,
\;\;\;\;\;\;\;\;\;\;
\omega_r=e_1+\cdots +e_r,
\end{equation}
where $P_\lambda (x)$ denotes the renormalized Askey-Wilson, Wilson, continuous
Hahn or Jacobi type polynomial of the form \eqref{reno}
and $\hat{E}_r(x)$ is the corresponding symmetric function multiplying
$P_\lambda (x)$
in the l.h.s. of the Recurrence relation \eqref{recr}.
(In all four cases of interest the functions $\hat{E}_r(x)$
are real for parameters in the domains of Section~\ref{sec2}
and we hence have \eqref{ident} trivially.)
If we work out both sides of \eqref{ident} by replacing
$\hat{E}_r P_\lambda$ and $\hat{E}_r P_{\lambda +\omega_r}$ by
the corresponding r.h.s. of \eqref{recr} and use the orthogonality
of the polynomials (Theorem~\ref{orthothm}),
then we arrive at the following relation between
$\langle P_\lambda ,P_\lambda \rangle_\Delta$ and
$\langle P_{\lambda +\omega_r} ,p_{\lambda +\omega_r} \rangle_\Delta$
\begin{eqnarray}\label{rel}
& &\hat{V}_{\{1,\ldots ,r\} ,\{r+1,\ldots ,n\} }(\rho +\lambda )
\langle P_{\lambda +\omega_r} ,P_{\lambda +\omega_r} \rangle_\Delta
= \\
& & \makebox[12em]{} \hat{V}_{\{1,\ldots ,r\} ,\{r+1,\ldots ,n\} }
(-\rho -\lambda -\omega_r)
\langle P_\lambda ,P_\lambda \rangle_\Delta . \nonumber
\end{eqnarray}
(Recall that $\hat{U}_{K,p}=1$ for $p=0$ and that
$\hat{V}_{\{1,\ldots ,r\} ,\{r+1,\ldots ,n\} }$ is taken to be
in accordance
with the definition in \eqref{recr} with
all signs $\varepsilon_j$, $j\in J$ being positive.)
In principle one can use this relation to obtain for each $\lambda\in\Lambda$
the squared norm of $P_\lambda(x)$ in terms of the squared norm of
$P_0(x)(=1)$ by writing $\lambda$ as a positive linear
combination of the (fundamental weight) vectors
$\omega_1,\ldots ,\omega_n$ and then apply
\eqref{rel} iteratively by walking to the weight $\lambda$ starting from
the zero weight $(0,\ldots ,0)$ through the successive addition of
fundamental weight vectors $\omega_r$.
Indeed, it is not difficult to verify that the combinatorial
structure of
the coefficients $\hat{V}_{\{1,\ldots ,r\} ,\{r+1,\ldots ,n\} }$ is
such that the result does not depend on the order in which
the fundamental weight vectors $\omega_r$ are added, i.e., the result
is independent of the chosen path from $(0,\ldots ,0)$ to $\lambda$.
This hinges on the easily inferred (combinatorial) identity
\begin{eqnarray*}
&&\frac{\hat{V}_{\{1,\ldots ,r\} ,\{r+1,\ldots ,n\} }(-x-\omega_s-\omega_r)}
     {\hat{V}_{\{1,\ldots ,r\} ,\{r+1,\ldots ,n\} }(x+\omega_s)}
\frac{\hat{V}_{\{1,\ldots ,s\} ,\{s+1,\ldots ,n\} }(-x-\omega_s)}
     {\hat{V}_{\{1,\ldots ,s\} ,\{s+1,\ldots ,n\} }(x)}
= \\ &&\makebox[3em]{}
\frac{\hat{V}_{\{1,\ldots ,s\} ,\{s+1,\ldots ,n\} }(-x-\omega_r-\omega_s)}
     {\hat{V}_{\{1,\ldots ,s\} ,\{s+1,\ldots ,n\} }(x+\omega_r)}
\frac{\hat{V}_{\{1,\ldots ,r\} ,\{r+1,\ldots ,n\} }(-x-\omega_r)}
     {\hat{V}_{\{1,\ldots ,r\} ,\{r+1,\ldots ,n\} }(x)} ,
\end{eqnarray*}
which expresses the fact that the result for
the quotient of $\langle P_{\lambda +\omega_r +\omega_s},
P_{\lambda+\omega_r +\omega_s} \rangle_\Delta $ and
$\langle P_\lambda ,P_\lambda\rangle_\Delta $
computed via \eqref{rel}
does not depend on the
order in which $\omega_r$ and $\omega_s$ are added
(as it clearly should not).

To write down the answer for
$\langle P_\lambda ,P_\lambda\rangle_\Delta /\langle 1,1\rangle_\Delta$
for general $\lambda\in\Lambda$ we introduce the functions
\begin{equation}\label{Deltahpm}
\hat{\Delta}_\pm (x) = \prod_{1\leq j< k \leq n}
\hat{d}_{v,\pm} (x_j+x_k)\, \hat{d}_{v,\pm}(x_j-x_k)
\prod_{1\leq j\leq n} \hat{d}_{w,\pm}(x_j),
\end{equation}
with $\hat{d}_{v,\pm} (z)$ and $\hat{d}_{w,\pm}(z)$  ($\neq 0$)
satisfying the difference equations
\begin{eqnarray}\label{diffv}
\hat{d}_{v,+} (z+1)=\hat{v}(z)\, \hat{d}_{v,+} (z),\;\;\;\;\;\;\;\;
\hat{d}_{v,-} (z+1)=\hat{v}(-z-1)\, \hat{d}_{v,-} (z), \\
\hat{d}_{w,+} (z+1)=\hat{w}(z)\, \hat{d}_{w,+} (z),\;\;\;\;\;\;\;\;
\hat{d}_{w,-} (z+1)=\hat{w}(-z-1)\, \hat{d}_{w,-} (z) ,\label{diffw}
\end{eqnarray}
where $\hat{v}(z)$ and $\hat{w}(z)$ are taken to be
the same as in Section~\ref{sec6}.
It is immediate from the difference equations \eqref{diffv}, \eqref{diffw}
that
\begin{eqnarray}\label{ado1}
\frac{\hat{\Delta}_+ (x+\omega_r)}{\hat{\Delta}_+(x)} &=&
\hat{V}_{ \{ 1,\ldots ,r\} ,\, \{ r+1,\ldots ,n\} } (x) ,\\
\label{ado2}
\frac{\hat{\Delta}_-(x+\omega_r)}{\hat{\Delta}_-(x)} &=&
\hat{V}_{ \{ 1,\ldots ,r\} ,\, \{ r+1,\ldots ,n\} }(-x-\omega_r ).
\end{eqnarray}
With the aid of the properties \eqref{ado1} and \eqref{ado2}
we can rewrite Relation~\eqref{rel} in the form
\begin{equation}\label{rewr}
\langle P_\lambda ,P_\lambda \rangle_\Delta \;
\frac{\hat{\Delta}_+ (\rho+\lambda)}{\hat{\Delta}_- (\rho+\lambda)}
=
\langle P_{\lambda +\omega_r} ,P_{\lambda +\omega_r} \rangle_\Delta \;
\frac{\hat{\Delta}_+ (\rho+\lambda +\omega_r)}
     {\hat{\Delta}_- (\rho+\lambda +\omega_r)} .
\end{equation}
By using the fact
that the fundamental weight vectors $\omega_1,\ldots ,\omega_n$
positively generate the integral cone $\Lambda $ \eqref{cone},
one deduces from this equation that the quotient
$\langle P_\lambda ,P_\lambda \rangle_\Delta\;
\hat{\Delta}_+ (\rho+\lambda)/\hat{\Delta}_- (\rho+\lambda)$ in the
l.h.s. of \eqref{rewr} does not depend on the choice of
$\lambda\in \Lambda$.
Comparing with its evaluation in $\lambda =0$ then entails
\begin{equation}\label{evalP}
\frac{\langle P_\lambda ,P_\lambda \rangle_\Delta}
     {\langle 1 , 1\rangle_\Delta }
=
\frac{\hat{\Delta}_- (\rho+\lambda )\: \hat{\Delta}_+ (\rho) }
     {\hat{\Delta}_+ (\rho+\lambda )\: \hat{\Delta}_- (\rho)}.
\end{equation}
As we will see below (and is suggested by the notation), it turns out
that the function $\Delta_+ (x)$ in the present section
coincides with that of Section~\ref{sec6}.
So, by combining \eqref{evalP} with \eqref{reno} we obtain
\begin{equation}\label{evalp}
\frac{\langle p_\lambda ,p_\lambda \rangle_\Delta}
     {\langle 1 , 1\rangle_\Delta } =
|c|^{-2|\lambda |} \;
\frac{\hat{\Delta}_+ (\rho+\lambda )\: \hat{\Delta}_- (\rho+\lambda ) }
     {\hat{\Delta}_+ (\rho )\: \hat{\Delta}_- (\rho ) }.
\end{equation}
We will now list for each of our four families AW, W, cH and J the associated
functions $\hat{d}_{v,\pm} (z)$ and $\hat{d}_{w,\pm} (z)$, and
formulate the corresponding
evaluation theorem for the quotient of
$\langle p_\lambda ,p_\lambda \rangle_\Delta $ and
$\langle 1 , 1\rangle_\Delta$. The proof of this theorem
boils in each case down to verifying that
$\hat{d}_{v,\pm} (z)$ and $\hat{d}_{w,\pm} (z)$ indeed satisfy
the difference equations \eqref{diffv} and \eqref{diffw}.

In  the case of multivariable Jacobi type polynomials
the value of the integral for the squared norms
$\langle p_\lambda ,p_\lambda
\rangle_\Delta$
was computed (in essence) by Opdam
\cite{opd:some} with the aid of shift operators
(see also \cite{hec-sch:harmonic}).
Formulas for the squared norms
of the Askey-Wilson type polynomials
were first conjectured by Macdonald \cite{mac:orthogonal,mac:some2} and
then proven by Cherednik \cite{che:double}
for special parameters
related to the reduced root systems (and admissible pairs of the
form $(R,R^\vee)$)
using a generalization of the shift operator approach.
Recently, Macdonald announced a further extension
of these methods to the case of general Askey-Wilson
parameters \cite{mac:affine}.

In \cite{die:self-dual} the present author combined a proof of \eqref{evalp}
using the recurrence relations
along the lines sketched above with Gustafson's constant
term formula \cite{gus:generalization,kad:proof}
for $\langle 1 ,1 \rangle_{\Delta^{AW}} $, which
led to an alternative derivation
of Macdonald's formula for the value
of $\langle p_\lambda^{AW} ,p_\lambda^{AW} \rangle_{\Delta^{AW}} $
(at present in the case of parameters satisfying Condition~\eqref{self-dualc}).
The (Askey-Wilson) case is included here just
for the sake of completeness. Also for the Jacobi type
polynomials our proof of
\eqref{evalp} with the aid of recurrence relations rather than shift
operators provides, when combined with
the previously derived expression for the constant term
$\langle 1 ,1 \rangle_{\Delta^{J}} $ \cite{mac:some1,sel:bemerkninger},
an alternative way
to demonstrate the validity of the known evaluation formula for
$\langle p_\lambda^{J} ,p_\lambda^{J} \rangle_{\Delta^{J}}$.

\subsection{Askey-Wilson type}
In the case of Askey-Wilson type polynomials
the functions
$\hat{\Delta}_\pm (x)=\hat{\Delta}_\pm^{AW} (x)$ of the form
in \eqref{Deltahpm} are
characterized by  $\hat{d}_{v,+}^{AW} (z)$ \eqref{dhAWvp},
$\hat{d}_{w,+}^{AW} (z)$ \eqref{dhAWwp} and
\begin{eqnarray}
\hat{d}^{AW}_{v,-}(z)& =&t^{z/2}
\frac{(q^{z+1}; q)_\infty}
     {( t^{-1} q^{z+1}; q)_\infty},    \\[1ex]
\hat{d}^{AW}_{w,-}(z)&=&
(\hat{t}_0 \hat{t}_1 \hat{t}_2 \hat{t}_3q^{-1})^{z/2}
\frac{(q^{2 z+1};\,
                                    q)_\infty    }
     {(\hat{t}^{-1}_0 q^{z+1},\,
       \hat{t}^{-1}_1 q^{z+1},\,
       \hat{t}^{-1}_2 q^{z+1},\,
       \hat{t}^{-1}_3 q^{z+1};\,
                                    q)_\infty    }
\end{eqnarray}
(with $\hat{t}_r$ given by \eqref{hattr}).
Formula \eqref{evalp} leads us to the following theorem for this case.
\begin{thm}[\cite{die:self-dual}]
Let us assume Askey-Wilson parameters with values taken
from the domain indicated in Section~\ref{AWdef}
such that the recurrence relations of Theorem~\ref{recrAWthm} hold.
Then one has
\begin{equation}\label{evalAWp}
\frac{\langle p^{AW}_\lambda ,p^{AW}_\lambda \rangle_{\Delta^{AW}}}
     {\langle 1 , 1\rangle_{\Delta^{AW}} } =
\frac{\hat{\Delta}^{AW}_+ (\rho^{AW}+\lambda )\:
      \hat{\Delta}^{AW}_- (\rho^{AW}+\lambda ) }
     {\hat{\Delta}^{AW}_+ (\rho^{AW} )\: \hat{\Delta}^{AW}_- (\rho^{AW} ) }
\end{equation}
(with $\rho^{AW}$ given by \eqref{rhoAW}).
\end{thm}
To complete the proof of the theorem it suffices to infer
that the functions $\hat{d}^{AW}_{v,\pm}(z)$ and $\hat{d}^{AW}_{w,\pm}(z)$
indeed satisfy the corresponding
difference equations of the form in \eqref{diffv} and \eqref{diffw};
a fact not difficult to deduce from the standard relation
for the $q$-shifted factorials
$(a;q)_\infty =(1-a)\, (aq;q)_\infty$.

For $n=1$ the r.h.s. of \eqref{evalAWp} reduces to
an expression of the form
$$\langle p^{AW}_l ,p^{AW}_l \rangle_{\Delta^{AW}}/
\langle p^{AW}_0 ,p^{AW}_0 \rangle_{\Delta^{AW}}$$ with
(cf. \cite{ask-wil:some,koe-swa:askey-scheme} and
recall our normalization in \eqref{bhrepAW})
\begin{equation}
\langle p^{AW}_l ,p^{AW}_l \rangle_{\Delta^{AW}}=
\frac{ 4\pi\alpha^{-1}\; (t_0t_1t_2t_3q^{2l};q)_\infty }
     { (t_0t_1t_2t_3q^{l-1};q)_l\: (q^{l+1};q)_\infty \:
       \prod_{0\leq r< s\leq 3} (t_rt_s q^l;q)_\infty} .
\end{equation}

\subsection{Wilson type}
The relevant functions
$\hat{\Delta}_\pm (x)=\hat{\Delta}_\pm^{W} (x)$ of the form
in \eqref{Deltahpm} are
determined by $\hat{d}_{v,+}^{W} (z)$,
$\hat{d}_{w,+}^{W} (z)$ \eqref{dhWvwp} and
\begin{equation}\label{dhWvwm}
\hat{d}^{W}_{v,-}(z)=
\frac{\Gamma (-\nu +z+1)}
     {\Gamma (z+1)},    \;\;\;\;\;\;\;\;\;\;\;\;\;\;
\hat{d}^{W}_{w,-}(z)=
\frac{\prod_{0\leq r\leq 3}\Gamma (-\hat{\nu}_r +z+1)   }
     {\Gamma (2z+1) }
\end{equation}
(with $\hat{\nu}_r$ given by \eqref{nuhr}). The difference equations
of the form in
\eqref{diffv} and \eqref{diffw} are easily verified using the
standard functional equation $\Gamma (z+1) =z \Gamma (z)$ for
the gamma function. The normalization theorem
becomes in the present situation.

\begin{thm}Let us assume Wilson parameters with values taken
from the domain indicated in Section~\ref{Wdef}
such that the recurrence relations of Theorem~\ref{recrWthm} hold.
Then one has
\begin{equation}\label{evalWp}
\frac{\langle p^{W}_\lambda ,p^{W}_\lambda \rangle_{\Delta^{W}}}
     {\langle 1 , 1\rangle_{\Delta^{W}} } =
\frac{\hat{\Delta}^{W}_+ (\rho^{W}+\lambda )\:
      \hat{\Delta}^{W}_- (\rho^{W}+\lambda ) }
     {\hat{\Delta}^{W}_+ (\rho^{W} )\: \hat{\Delta}^{W}_- (\rho^{W} ) }
\end{equation}
(with $\rho^W$ given by \eqref{rhoW}).
\end{thm}
For $n=1$ the r.h.s. of \eqref{evalWp}
reduces to
$$\langle p^{W}_l ,p^{W}_l \rangle_{\Delta^{W}}/
\langle p^{W}_0 ,p^{W}_0 \rangle_{\Delta^{W}}$$
with (cf. \cite{koe-swa:askey-scheme} and
recall our normalization in \eqref{hrepW})
\begin{equation}
\langle p^{W}_l ,p^{W}_l \rangle_{\Delta^{W}}=
\frac{4\pi \; l! \:\prod_{0\leq r< s\leq 3} \Gamma (\nu_r +\nu_s +l) }
{  (\nu_0+\nu_1+\nu_2+\nu_3+l-1)_l\: \Gamma (\nu_0+\nu_1+\nu_2+\nu_3+2l)}.
\end{equation}

\subsection{Continuous Hahn type}
The functions
$\hat{\Delta}_\pm (x)=\hat{\Delta}_\pm^{cH} (x)$ of the form
in \eqref{Deltahpm} are now
determined by $\hat{d}_{v,+}^{cH} (z)$,
$\hat{d}_{w,+}^{cH} (z)$ \eqref{dhcHvwp} and
\begin{equation}\label{dhcHvwm}
\hat{d}^{cH}_{v,-}(z)=
\frac{\Gamma (-\nu +z+1)}
     {\Gamma (z+1)},    \;\;\;\;\;\;\;\;\;\;\;\;\;\;
\hat{d}^{cH}_{w,-}(z)=
\frac{\prod_{0\leq r\leq 2}\Gamma (-\hat{\nu}_r +z+1)   }
     {\Gamma (2z+1) }
\end{equation}
(with $\hat{\nu}_r$ taken from \eqref{nuhrcH}).
The difference equations \eqref{diffv}, \eqref{diffw} follow again
from the standard difference equation for the gamma function.
The normalization theorem reads in this case.

\begin{thm}Let us assume continuous Hahn parameters with values taken
from the domain indicated in Section~\ref{cHdef}
such that the recurrence relations of Theorem~\ref{recrcHthm} hold.
Then one has
\begin{equation}\label{evalcHp}
\frac{\langle p^{cH}_\lambda ,p^{cH}_\lambda \rangle_{\Delta^{cH}}}
     {\langle 1 , 1\rangle_{\Delta^{cH}} } =
\frac{\hat{\Delta}^{cH}_+ (\rho^{cH}+\lambda )\:
      \hat{\Delta}^{cH}_- (\rho^{cH}+\lambda ) }
     {\hat{\Delta}^{cH}_+ (\rho^{cH} )\: \hat{\Delta}^{cH}_- (\rho^{cH} ) }
\end{equation}
(with $\rho^{cH}$ given by \eqref{rhocH}).
\end{thm}
For $n=1$ the r.h.s. of \eqref{evalcHp} becomes
$$   \langle p^{cH}_l ,p^{cH}_l \rangle_{\Delta^{cH}}/
\langle p^{cH}_0 ,p^{cH}_0 \rangle_{\Delta^{cH}} $$ with
(cf. \cite{koe-swa:askey-scheme} and recall our
normalization in \eqref{hrepcH})
\begin{eqnarray}
&& \langle p^{cH}_l ,p^{cH}_l \rangle_{\Delta^{cH}}=\\
&&\makebox[4em]{}
\frac{2\pi \; l! \:\prod_{r,s\in \{ 0, 1\} } \Gamma (\nu_r^+ +\nu_s^- +l) }
     {  (\nu_0^+ +\nu_0^- +\nu_1^+ +\nu_1^- +l-1)_l\:
      \Gamma (\nu_0^+ +\nu_0^- +\nu_1^+ +\nu_1^- +2l)}.\nonumber
\end{eqnarray}

\subsection{Jacobi type}
For the Jacobi type polynomials the functions
$\hat{\Delta}_\pm (x)=\hat{\Delta}_\pm^{J} (x)$ of the form
in \eqref{Deltahpm} are
determined by $\hat{d}_{v,+}^{J} (z)$,
$\hat{d}_{w,+}^{J} (z)$ \eqref{dhJvwp} and
\begin{equation}\label{dhJvwm}
\hat{d}^{J}_{v,-}(z)=
\frac{\Gamma (-\nu +z+1)}
     {\Gamma (z+1)},    \;\;\;\;\;\;\;\;\;\;
\hat{d}^{J}_{w,-}(z)=
\frac{\Gamma (-\hat{\nu}_0 +z+1) \Gamma (-\hat{\nu}_1 +z+1)  }
     {\Gamma (2z+1) }  .
\end{equation}
(with $\hat{\nu}_0$, $\hat{\nu}_1$  taken from \eqref{nuhrJ}).
The difference equations \eqref{diffv}, \eqref{diffw} follow again
from the difference equation for the gamma function.
The normalization theorem reads in this case.

\begin{thm}Let us assume Jacobi parameters with values taken
from the domain indicated in Section~\ref{Jdef}.
Then one has
\begin{equation}\label{evalJp}
\frac{\langle p^{J}_\lambda ,p^{J}_\lambda \rangle_{\Delta^{J}}}
     {\langle 1 , 1\rangle_{\Delta^{J}} } = 2^{4|\lambda |}\;
\frac{\hat{\Delta}^{J}_+ (\rho^{J}+\lambda )\:
      \hat{\Delta}^{J}_- (\rho^{J}+\lambda ) }
     {\hat{\Delta}^{J}_+ (\rho^{J} )\: \hat{\Delta}^{J}_- (\rho^{J} ) }
\end{equation}
(with $\rho^{J}$ given by \eqref{rhoJ}).
\end{thm}
For $n=1$ the r.h.s. of \eqref{evalJp}
becomes
$$ \langle p^{J}_l ,p^{J}_l \rangle_{\Delta^{J}}/
\langle p^{J}_0 ,p^{J}_0 \rangle_{\Delta^{J}}$$
with
(cf. \cite{abr-ste:handbook,koe-swa:askey-scheme} and recall our
normalization in \eqref{hrepJ})
\begin{equation}
\langle p^{J}_l ,p^{J}_l \rangle_{\Delta^{J}}=
\frac{2^{4l+1}\alpha^{-1}\,
       l!\; \Gamma (\nu_0+1/2) \Gamma (\nu_1+1/2)          }
     {  (\nu_0+\nu_1 +l)_l\:     \Gamma (\nu_0+\nu_1 +2l+1) } .
\end{equation}

\newpage
\section*{Acknowledgments}
The author would like to thank Prof. T. Oshima for the kind
hospitality at the University of Tokyo.
Thanks are also due to J. V. Stokman for pointing out the
usefulness of Formula~\eqref{urep} in studying the transition
to degenerate cases of the multivariable
Askey-Wilson polynomials.

\bibliographystyle{amsplain}

\begin{thebibliography}{0000}

\bibitem[AS]{abr-ste:handbook} M. Abramowitz and I. A. Stegun (eds.),
{\em Handbook of mathematical functions}, Dover Publications,
New York, 1972 (9th printing).

\bibitem[AW]{ask-wil:some} R. Askey and J. Wilson, {\em Some basic
hypergeometric orthogonal polynomials that generalize Jacobi polynomials},
Mem. Amer. Math. Soc. {\bf 54} (1985), no. 319.

\bibitem[BO]{bee-opd:certain} R. J. Beerends and E. M. Opdam, {\em
Certain hypergeometric series related to the root system $BC$},
Trans. Amer. Math. Soc. {\bf 339} (1993), 581-609.

\bibitem[C1]{che:double} I. Cherednik, {\em Double affine Hecke algebras and
Macdonald's conjectures}, Ann. Math. {\bf 141} (1995), 191-216.

\bibitem[C2]{che:macdonalds} \bysame, {\em Macdonald's evaluation
conjectures and difference Fourier transform}, Invent. Math.
{\bf 122} (1995), 119-145.

\bibitem[De]{deb:systeme} A. Debiard, {\em Syst\`eme diff\'erentiel
hyperg\'eom\'etrique et parties radiales des op\'erateurs
invariants des espaces sym\'etriques de type $BC_p$}, in:
S{\'e}minaire d'Alg{\`e}bre Paul Dubreil
et Marie-Paule Malliavin (M.-P. Malliavin, ed.), Lecture Notes in Math.,
vol. 1296, Springer, Berlin, 1988, pp. 42-124.

\bibitem[D1]{die:commuting} J. F. van Diejen, {\em Commuting difference
operators with polynomial eigenfunctions}, Compositio Math. {\bf 95}
(1995), 183-233.

\bibitem[D2]{die:difference} \bysame, {\em Difference Calogero-Moser systems
and finite Toda chains}, J. Math. Phys. {\bf 36} (1995), 1299-1323.

\bibitem[D3]{die:multivariable} \bysame, {\em Multivariable continuous Hahn and
Wilson polynomials related to integrable difference systems},
J. Phys. A: Math. Gen. {\bf 28} (1995), L369-L374.

\bibitem[D4]{die:diagonalization} \bysame, {\em On the diagonalization of
difference Calogero-Sutherland systems}, in:
Proceedings of the workshop on symmetries and
integrability of difference equations (D. Levi, L. Vinet, and P. Winternitz,
eds.), CRM Proceedings and Lecture Notes (to appear).

\bibitem[D5]{die:self-dual} \bysame, {\em Self-dual
Koornwinder-Macdonald polynomials}, Invent. Math. (to appear).

\bibitem[Du]{dun:differential} C. F. Dunkl, {\em Differential-difference
operators associated to reflection groups},
Trans. Amer. Math. Soc. {\bf 311} (1989), 167-183.

\bibitem[GR]{gas-ram:basic} G. Gasper and M. Rahman, {\em Basic
hypergeometric series}, Cambridge University Press, Cambridge, 1990.

\bibitem[Gu]{gus:generalization} R. A. Gustafson, {\em A generalization of
Selberg's beta integral}, Bull. Amer. Math. Soc. (N.S.) {\bf 22} (1990),
97-105.

\bibitem[H]{hec:elementary} G. J. Heckman, {\em An elementary
approach to the hypergeometric shift operator of Opdam},
Invent. Math. {\bf 103} (1991), 341-350.

\bibitem[HS]{hec-sch:harmonic} G. J. Heckman and H.
Schlichtkrull, {\em Harmonic analysis and special functions
on symmetric spaces}, Perspectives in Math., vol. 16,
Academic Press, San Diego, 1994.

\bibitem[Ka]{kad:proof} K. W. J. Kadell, {\em A proof of the
$q$-Macdonald-Morris conjecture for $BC_n$},
Mem. Amer. Math. Soc. {\bf 108} (1994), no. 516.

\bibitem[KS]{koe-swa:askey-scheme} R. Koekoek and R. F. Swarttouw,
{\em The Askey-scheme of hypergeometric orthogonal polynomials and its
$q$-analogue}, Math. report Delft Univ. of Technology 94-05, 1994.


\bibitem[K]{koo:askey} T. H. Koornwinder,  {\em Askey-Wilson
polynomials for root systems of type BC}, in:
Hypergeometric functions on domains of positivity, Jack polynomials, and
applications (D. St. P. Richards, ed.),
Contemp. Math., vol. 138, Amer. Math. Soc., Providence, R. I., 1992,
pp. 189-204.

\bibitem[M1]{mac:some1} I. G. Macdonald, {\em Some conjectures for root
systems}, SIAM J. Math. Anal. {\bf 13}
(1982), 988-1007.

\bibitem[M2]{mac:orthogonal} \bysame, {\em Orthogonal polynomials
associated with root systems},  unpublished manu\-script, 1988.

\bibitem[M3]{mac:some2} \bysame, {\em Some conjectures for Koornwinder's
orthogonal polynomials}, unpublished manu\-script, 1991.

\bibitem[M4]{mac:symmetric} \bysame, {\em Symmetric functions and Hall
polynomials}, 2nd edition, Clarendon Press, Oxford, 1995.

\bibitem[M5]{mac:affine} \bysame, {\em Affine Hecke algebras and orthogonal
polynomials}, S\'eminaire Bourbaki {\bf 47} (1995), no. 797, 1-18.

\bibitem[OOS]{och-osh-sek:commuting} H. Ochiai, T. Oshima, and H. Sekiguchi,
{\em Commuting families of symmetric differential operators},
Proc. Japan Acad. Ser. A Math. Sci. {\bf 70} (1994), 62-66.

\bibitem[OP]{ols-per:quantum} M. A. Olshanetsky and A. M. Perelomov,
{\em Quantum integrable systems related to Lie algebras},
Phys. Reps. {\bf 94} (1983), 313-404.

\bibitem[Op]{opd:some} E. M. Opdam, {\em Some applications of hypergeometric
shift operators}, Invent. Math. {\bf 98} (1989), 1-18.

\bibitem[OS]{osh-sek:commuting}  T. Oshima and H. Sekiguchi,
{\em Commuting families of differential operators invariant under the action
of a Weyl group}, J. Math. Sci. Univ. Tokyo {\bf 2} (1995), 1-75.

\bibitem[R1]{rui:complete} S. N. M. Ruijsenaars, {\em Complete integrability of
relativistic Calogero-Moser systems and elliptic function identities},
Commun. Math. Phys. {\bf 110} (1987), 191-213.

\bibitem[R2]{rui:finite} \bysame, {\em Finite-dimensional
soliton systems}. in.: Integrable and superintegrable systems
(B. Kupershmidt, ed.), World Scientific, Singapore, 1990,
pp. 165-206.

\bibitem[Se]{sel:bemerkninger} A. Selberg, {\em Bemerkninger om et
multipelt integral}, Norsk Mat. Tidsskr. {\bf 26} (1944), 71-78
(Collected papers, vol. 1, Springer, Berlin, 1989, pp. 204-213).

\bibitem[S]{sto:multivariable} J. V. Stokman, {\em Multivariable big
and little $q$-Jacobi polynomials}, Math. preprint Univ. of Amsterdam
95-16, 1995.

\bibitem[SK]{sto-koo:limit} J. V. Stokman and T. H. Koornwinder, {\em Limit
transitions for $BC$ type multivariable orthogonal polynomials},
Math. preprint Univ. of Amsterdam 95-19, 1995.

\bibitem[V]{vre:formulas} L. Vretare, {\em Formulas for elementary spherical
functions and generalized Jacobi polynomials}, SIAM J. Math. Anal. {\bf 15}
(1984), 805-833.

\end{thebibliography}

\end{document}